# On the constructions of the skew Brownian motion

## Antoine Lejay[*]

*IECN,*
*Campus scientifique,*
*BP 239,*
*54506 Vandœuvre-lès-Nancy* CEDEX, *France*
*e-mail:* `Antoine.Lejay@iecn.u-nancy.fr`

**Abstract:** This article summarizes the various ways one may use to construct the Skew Brownian motion, and shows their connections. Recent applications of this process in modelling and numerical simulation motivates this survey. This article ends with a brief account of related results, extensions and applications of the Skew Brownian motion.

## 1. Introduction

The Skew Brownian motion appeared in the '70 in [44, 87] as a natural generalization of the Brownian motion: it is a process that behaves like a Brownian motion except that the sign of each excursion is chosen using an independent Bernoulli random variable of parameter $p$. As a sequel, if $X$ is a Skew Brownian motion of parameter $p$ (that is, $p$ is the probability that the excursion is positive), then $\mathbb{P}[X_t > 0] = p$ for any $t > 0$, while outside 0, $X_t$ behaves like a Brownian motion. As shown in [41], this process is a semi-martingale which is a strong solution to some Stochastic Differential Equation (SDE) with local time

$$X_t = x + B_t + (2p - 1)L_t^0(X),$$

where $L_t^0(X)$ is the local time at 0 of $X$.

At the same time, N. Portenko constructed in [69, 70] a process whose infinitesimal generator has a singular drift concentrated on some hypersurface. This is a way to model some permeable barrier. Indeed, he considers diffusion

[*]This research was partially funded by the French Groupe de Recherche MOMAS (ANDRA, BRGM, CEA, CNRS and EDF).

whose infinitesimal generator is

$$L = \sum_{i,j=1}^{d} \frac{1}{2} a_{i,j} \frac{\partial^2}{\partial x_i \partial x_j} + \sum_{i=1}^{d} b_i \frac{\partial}{\partial x_i} + \sum_{i=1}^{d} q \delta_S N_i(x) \frac{\partial}{\partial x_i}$$

where $N = (N_1, \ldots, N_d)$ is the vector conormal to a smooth hyper-surface $S$, $a$, $b$ and $q$ are continuous and $|q(x)| < 1$ for all $x \in S$. He showed that the solution to $\frac{\partial u}{\partial t} = Lu$ with $u(0, x) = \varphi(x)$ is also solution to the PDE

$$\begin{cases} \dfrac{\partial u(t,x)}{\partial t} = Lu(t,x) \text{ on } \mathbb{R}_+^* \times \mathbb{R}^d \setminus S, \\ (1 + q(x))N_+(x) \cdot \nabla u(t,x) = (1 - q(x))N_-(x) \cdot \nabla u(t,x) \\ u(t,\cdot) \text{ is continuous on } S, \\ u(0,x) = \varphi(x) \end{cases}$$

where $N_+$ and $N_-$ are the inner and outer conormal to $S$. The condition on the left and right flux on the surface $S$ is called a *transmission condition* or a *flux condition*.

In dimension one with $a = 1$ and $b = 0$, this process is the Skew Brownian motion. Thus, the infinitesimal generator of this diffusion process has a rather natural interpretation, since it corresponds to half the Laplace operator plus a "generalized drift" given by a Dirac mass at 0 with a coefficient that corresponds to its skewness. Besides, the effect of the drift is translated into the PDE as a transmission condition. It was also noted in [92] that SDEs with local time may be used to model diffusion processes in a media with permeable barrier.

Since, a few works have used the Skew Brownian motion as a tool for solving applied problems: in astrophysics [98], in ecology [15], in homogenization [50, 92], in geophysics [56, 57] and more recently in finance [19, 20, 21, 22],...

A recent series of works [27, 56, 58, 61] shows that the properties of the Skew Brownian motion may be used in a systematic way to provide different schemes to simulate diffusion processes generated by

$$L = \frac{\rho}{2} \frac{\mathrm{d}}{\mathrm{d}x} \left( a \frac{\mathrm{d}}{\mathrm{d}x} \right) + b \frac{\mathrm{d}}{\mathrm{d}x}, \tag{1}$$

with $a$ and $\rho$ have discontinuous points of first kind. All these Monte Carlo methods rely on a probabilistic interpretation of the transmission conditions at the points where the coefficients are discontinuous and appropriate changes of scales.

Indeed, the diffusion process $X$ generated by $L$ is solution to an SDE with local time of type

$$\mathrm{d}X_t = \sigma(X_t)\,\mathrm{d}B_t + b(X_t)\,\mathrm{d}t + \int_{\mathbb{R}} \nu(\,\mathrm{d}x)\,\mathrm{d}L_t^x(X), \tag{2}$$

where $\nu$ is a finite, signed measure that has a mass at the points where $a$ or $\rho$ are discontinuous.

SDEs of type (2) have been studied first in [24, 51, 52] (see also [9, 25] for some extensions), where existence and uniqueness of a strong solution is proved under some rather general conditions ($\sigma$ shall be uniformly elliptic, bounded and of finite variation, $b$ shall be bounded, and $\nu$ shall be of finite mass with $|\nu(\{x\})| < 1$ for any point $x \in \mathbb{R}$).

This type of SDE generalized the usual SDEs and were studied with hypotheses on the coefficients that are as weak as possible: [9, 10, 30], ... In addition, it is worth noting that the articles [24, 52] also give an account on the "critical cases" between strong and weak solutions of SDEs, in the sense they deal with the weakest possible conditions on the coefficients to ensure the existence of a strong solution (see also [6] for an example of application to the theory of SDEs). In addition, some Dirichlet processes can be constructed as solutions of equations of this type [8, 25].

Another natural extension of the Skew Brownian motion is the *Walsh Brownian motion*. It is a diffusion process that moves on rays emanating from a single point and was introduced by J. Walsh in [87]. This diffusion process can be used locally as a description of a *spider martingale* [94, Sect. 17, p. 103], and is a special case of a diffusion on a graph, since we recover a transmission condition at each vertex. Diffusions on graphs are particularly important in the study of dynamical Hamiltonian systems as shown in the pioneering work of M. Freidlin and M. Weber [34], but could be useful for modelling many physical or biological diffusion phenomena. It is then possible to use Monte Carlo methods to solve that type of problems.

One has also to note that Walsh Brownian motion was useful to provide a counter-example to a natural question on Brownian filtration, as B. Tsirelson shown it in [85].

The goal of this article is then to summarize the different ways to construct a Skew Brownian motion (using PDEs, Dirichlet forms, approximations by diffusions with smooth coefficients and by random walks, scale functions and speed measures, SDE, excursions theory...) and their relationships. The last section presents quickly some extensions of the Skew Brownian motion and their applications to various fields.

### Notations

Classically, we denote by $\mathbb{R}_+$ (resp. $\mathbb{R}_-$) the set of non-negative (resp. non-positive) real numbers, and by $\mathbb{R}_+^*$ (resp. $\mathbb{R}_-^*$) the set $\mathbb{R}_+ \setminus \{0\}$ (resp. $\mathbb{R}_- \setminus \{0\}$).

The space of continuous functions from X to Y is denoted by $\mathcal{C}(X, Y)$.

The space of square integrable function $f$ on X is denoted by $L^2(X)$.

On $\mathbb{R}^d$ with $d \geq 1$, we denote by $H^1(\mathbb{R}^d)$ the completion of the space of smooth functions with compact support with respect to the norm

$$\|f\|_{H^1} = \sqrt{\int_{\mathbb{R}^d} |f(x)|^2 \, \mathrm{d}x + \sum_{i=1}^d \int_{\mathbb{R}^d} |\partial_{x_i} f(x)|^2 \, \mathrm{d}x}.$$

A function in $\mathrm{H}^1(\mathbb{R}^d)$ has then a square integrable generalized derivative $\partial_{x_i} f$ with respect to any coordinate $x_i$ in $\mathbb{R}^d$. Let us recall that when $d = 1$, for any function $f$ in $\mathrm{H}^1(\mathbb{R})$, there exists a continuous function $\tilde{f}$ such that $f(x) = \tilde{f}(x)$ almost everywhere.

The space $\mathrm{H}^2(\mathbb{R})$ contains the function $f$ in $\mathrm{H}^1(\mathbb{R})$ such that $\nabla f$ also belongs to $\mathrm{H}^1(\mathbb{R})$.

## 2. Differential operators with generalized coefficients

Let $\delta$ denote the Dirac function at 0. For the first time, we are interested in solving the parabolic PDE for some $q \in [-1, 1]$,

$$\begin{cases} \dfrac{\partial u(t,x)}{\partial t} = \frac{1}{2}\triangle u(t,x) + q\delta\nabla u(t,x), \\ u(0,x) = \varphi(x), \end{cases} \tag{3}$$

which is equivalent to construct the semi-group associated to the operator with a singular first-order differential term

$$L = \frac{1}{2}\triangle + q\delta\nabla. \tag{4}$$

In [69, 70], N. Portenko constructed the semi-group generated by $L$ and showed it is a Feller semi-group. In fact, N. Portenko works in $\mathbb{R}^m$ instead of $\mathbb{R}$, where $\delta$ is a Dirac mass of a surface $S$ smooth enough. Here, we restrict ourselves to the simpler form of $L$ given by (4), and we give a short account of the whole construction in Section 11.10.1 below.

### 2.1. Perturbation of the heat kernel

One wants to construct the semi-group $(Q_t)_{t>0}$, if it exists, generated by $L$. Let us denote by $(P_t)_{t>0}$ the semi-group generated by $\frac{1}{2}\triangle$ over $\mathbb{R}$. It is well known that $P_t$ has a density $p(t, x, y)$ given by the heat (or Gaussian) kernel: for any continuous, bounded function $f$,

$$P_t\varphi(x) = \int_{\mathbb{R}} p(t,x,y)\varphi(y)\,\mathrm{d}y \text{ with } p(t,x,y) = \frac{1}{\sqrt{2\pi t}}\exp\left(-\frac{|x-y|^2}{2t}\right).$$

It is natural to construct $Q_t$ as a perturbation of $P_t$, that is

$$Q_t\varphi(x) = P_t\varphi(x) + q\int_0^t p(t-\tau,x,0)\nabla Q_\tau\varphi(0)\,\mathrm{d}\tau. \tag{5}$$

The problem here is to give a meaning to $\nabla Q_t\varphi(0)$, since we will see that $\nabla Q_t\varphi$ is discontinuous at 0. However, let us remark that for any continuous function $f$ and any $\tau > 0$, $x \mapsto p(\tau, x, 0)f(0)$ is of class $\mathcal{C}^1(\mathbb{R}; \mathbb{R})$ and

$$\frac{\partial}{\partial x}\, p(\tau,x,0)f(0)|_{x=0} = 0. \tag{6}$$

So, if we inject in the right-hand side of (5) the value of $Q_t\varphi(x) \stackrel{\text{def.}}{=} u(t,x)$ given by (5), one gets with (6) that

$$Q_t\varphi(x) = u(t,x) = P_t\varphi(x) + q \int_0^t p(t-\tau,x,0)\nabla P_\tau\varphi(0)\,\mathrm{d}\tau. \tag{7}$$

We use (7) as the definition of $Q_t$.

### 2.2. Single layer potential

Another idea to construct the semi-group is to consider that the PDE (3) corresponds to the potential generated by a charge whose value is

$$V(t,\varphi) = q\frac{\nabla u(t,0+) + \nabla u(t,0-)}{2}$$

at the point 0. From standard results (See [31, 38] for example), one knows that the solution $Q_t\varphi(x) \stackrel{\text{def.}}{=} u(t,x)$ of such a problem is given on $\mathbb{R}_+^*$ and $\mathbb{R}_-^*$ by

$$Q_t\varphi(x) = P_t\varphi(x) + \int_0^t p(t-\tau,x,0)V(\tau,\varphi)\,\mathrm{d}\tau.$$

Yet $V(\tau,\varphi)$ also involves the value of $u$. Indeed, a computation similar to the one done previously leads also to (7), and thus

$$V(t,\varphi) = \nabla P_t\varphi(0). \tag{8}$$

The advantage of this formulation is that one also knows from standard results on potential theory (see [38] for example) that

$$\nabla Q_t\varphi(0\pm) = \int_\mathbb{R} \nabla p(t,0,y)\varphi(y)\,\mathrm{d}y \mp qV(\tau,\varphi)$$

which means that with (8),

$$\nabla Q_t\varphi(0+) = (1-q)\nabla P_t\varphi(0) \text{ and } \nabla Q_t\varphi(0-) = (1+q)\nabla P_t\varphi(0). \tag{9}$$

Then it is immediate that $(t,x) \mapsto Q_t\varphi(x)$ is of class $\mathcal{C}^{1,2}$ on $\mathbb{R}_+^* \times \mathbb{R}^*$, is continuous on $\mathbb{R}^* \times \mathbb{R}$, but $x \mapsto \nabla Q_t\varphi(x)$ is discontinuous at $x=0$. From (9), one gets

$$\alpha\nabla Q_t\varphi(0+) = (1-\alpha)\nabla Q_t\varphi(0-) \tag{10}$$

with

$$\alpha = \frac{1+q}{2}. \tag{11}$$

Hence, one can give a proper meaning to (3), since $(t,x) \mapsto Q_t\varphi(x)$ is a continuous solution to

$$\begin{cases} (t,x) \mapsto u(t,x) \in \mathcal{C}^{1,2}(\mathbb{R}_+^* \times \mathbb{R}^*; \mathbb{R}) \cap \mathcal{C}(\mathbb{R}_+^* \times \mathbb{R}; \mathbb{R}), \\ \dfrac{\partial u(t,x)}{\partial t} = \frac{1}{2}\triangle u(t,x) \text{ for any } (t,x) \in \mathbb{R}_+^* \times \mathbb{R}^*, \\ \alpha\nabla u(t,0+) = (1-\alpha)\nabla u(t,0-), \\ u(t,0+) = u(t,0-), \\ u(0,x) = \varphi(x). \end{cases} \tag{12}$$

**Proposition 1.** *When $\varphi$ is continuous and bounded, and $\alpha \in [0, 1]$, there exists a unique solution to (12) for which the maximum principle holds.*

*Proof.* The existence and uniqueness of the solution of (12) may be proved using the equivalence between this problem and the same PDE written in a variational form: see Section 3.1.

But we could also prove this result without using the notion of weak solutions. The existence of $u$ is proved, since we have constructed such a solution with the help of the semi-group $(Q_t)_{t>0}$.

Let us note also that the uniqueness of $u$ follows from the maximum principle. If the initial condition $\varphi$ is non-negative (resp. non-positive), then $u$ is non-negative (resp. non-positive). Hence, if $\varphi = 0$, then $u$ is both non-negative and non-positive and is then equal to 0. As (12) is linear, this yields the uniqueness of its solution.

If $\alpha = 0$ (resp. $\alpha = 1$), then (12) is equivalent to two heat equations, one with a Neumann boundary condition on $\mathbb{R}_-^* \times \{0\}$ (resp. $\mathbb{R}_+^* \times \{0\}$), and the other with a lateral Dirichlet boundary on $\mathbb{R}_+^* \times \{0\}$ (resp. $\mathbb{R}_-^* \times \{0\}$) which is specified, thanks to the continuity of $u$, to the value of $u$ on $\mathbb{R}_-^* \times \{0\}$ (resp. $\mathbb{R}_+^* \times \{0\}$). Hence, the maximum principle holds as it holds for both equations.

We deal now with the case $\alpha \in (0, 1)$.

Assume that $\varphi \leq 0$ and $\varphi$ has a compact support. Indeed, if $u$ is a solution of (12), then $u$ is solution to $\frac{\partial u}{\partial t}(t, x) = \frac{1}{2}\triangle u(t, x)$ both on $\mathbb{R}_+^* \times \mathbb{R}_+$ and $\mathbb{R}_+^* \times \mathbb{R}_-$, with the boundary condition $\varphi$ on $\mathbb{R}_+$ and $\mathbb{R}_-$, and the lateral boundary condition $u(t, 0)$ on $\mathbb{R}_+^* \times \{0\}$. The solution $u$ for each of these equations satisfies $u(t, x) \to 0$ when $|x| \to +\infty$ or $t \to +\infty$.

Thus, if $u > 0$ inside $\mathbb{R}_+^* \times \mathbb{R}$ (resp. $\mathbb{R}_-^* \times \mathbb{R}$) the maximum of $u$ is attained at a point $(t_0^+, x_0^+)$ of $\mathbb{R}_+^* \times \mathbb{R}_+$ (resp. $(t_0^-, x_0^-)$ of $\mathbb{R}_+^* \times \mathbb{R}_-$). If $x_0^\pm \neq 0$, then it follows from Theorem 2 in [38, Sect. 1, p. 38] that $u(t_0^\pm, 0) = u(t_0^\pm, x_0^\pm)$. In other words, the maximum of $u$ is then reached on the line $\mathbb{R}_+^* \times \{0\}$.

Let $t_0$ be a time such that $u(t_0, 0)$ is maximum. As $\alpha \nabla u(t_0, 0+) = (1 - \alpha)\nabla u(t_0, 0-)$, $\nabla u(t_0, 0+)$ and $\nabla u(t_0, 0-)$ have the same sign. If $\nabla u(t_0, 0+) > 0$ (resp. $\nabla u(t_0, 0-) < 0$), a Taylor development in $x$ around 0 shows that there exists some $\varepsilon > 0$ for which $u(t_0, x) > u(t_0, 0)$ for $x \in (0, \varepsilon)$ (resp. $x \in (-\varepsilon, 0)$, which contradicts that $u(t_0, 0)$ is a maximum of $u$. Yet Theorem 14 in [38, Sect. 1, p. 49] asserts that $\nabla u(t_0, 0\pm) \neq 0$. Hence, $u$ has no positive maximum in $\mathbb{R}_+^* \times \mathbb{R}$ and thus $u \leq 0$.

Of course, if $\varphi \geq 0$, then $-u$ is solution to (12) with $-\varphi$ as an initial condition and $u \geq 0$.

If $\varphi \geq 0$, then $-\varphi \leq 0$ and then $u \geq 0$, since $-u$ is solution to (12) with $-\varphi$ as an initial condition. To conclude, let us note that if $\varphi(x) = C$ for all $x \in \mathbb{R}$, then $u(t, x) = C$.

It remains to drop the assumption that $\varphi$ has a compact support. From (7), if $\varphi$ is continuous, $|\varphi|$ is bounded and $\varphi(x) = 0$ on $(-R, R)$, then $\sup_{(t,x) \in \mathbb{R}_+^* \times \mathbb{R}} |Q_t \varphi(x)|$ converges to 0 as $R \to \infty$. By combining this result with the previous one when $\varphi$ has a compact support, one easily gets that if $\varphi \leq C$ (resp. $\varphi \geq C$), then the solution $Q_t \varphi(x)$ of (12) satisfies $Q_t \varphi(x) \leq C$ (resp. $Q_t \varphi(x) \geq C$) for all

$(t, x) \in \mathbb{R}_+ \times \mathbb{R}$. □

The next proposition follows from the previous facts.

**Proposition 2** ([69]). *If $|q| \leq 1$ (or equivalently $\alpha \in [0, 1]$), then $(Q_t)_{t>0}$ is a Feller semi-group, where $Q_t$ is defined by (7) for any function $\varphi$ which is continuous and bounded on $\mathbb{R}$.*

### 2.3. Construction of a Skew Brownian motion

The consequences of the construction given by N. Portenko are summarized in the following theorem.

**Theorem 1.** *If $|q| \leq 1$, then $(Q_t)_{t>0}$ is the semi-group of a strong Markov process $(X_t, \mathcal{F}_t, \mathbb{P}_x; t \geq 0, x \in \mathbb{R})$ on a probability space $(\Omega, \mathcal{F}, \mathbb{P})$, which we call a Skew Brownian motion of parameter $\alpha$ (abbreviated by SBM($\alpha$)) with $\alpha = (1+q)/2$. This process is continuous and conservative. Besides, there exists a $(\mathcal{F}_t, \mathbb{P})$-Brownian motion $B$ and a continuous additive functional $\eta$ such that*

$$X_t = x + B_t + \eta_t, \ \forall t \geq 0, \ \mathbb{P}_x\text{-a.s..} \tag{13}$$

*The additive functional $\eta$ is of finite variation, and it variation increases only when the process is at $0$: $\int_0^t |\,\mathrm{d}\eta_s| = \int_0^t \mathbf{1}_{\{X_s=0\}} |\,\mathrm{d}\eta_s| \ \mathbb{P}_x\text{-a.s..}$*

Of course, it remains to prove that the infinitesimal generator of $X$ is $L$ given by (4). It is clear from (7) that $Q_t$ has a density $q(t, x, y)$ given by

$$q(t, x, y) = p(t, x, y) + \frac{1}{\pi} \int_0^t \frac{y}{\sqrt{\tau}\sqrt{t-\tau}} \exp\left(-\frac{x^2}{t-\tau} - \frac{y^2}{\tau}\right) \mathrm{d}\tau. \tag{14}$$

In fact, a more convenient expression of $q(t, x, y)$ will be given later in (17). For that, it could be shown that for any function $\psi$ continuous and bounded,

$$\lim_{t\to 0} \int_{\mathbb{R}} \psi(x) \left( \frac{1}{t} \int_{\mathbb{R}} q(t, x, y)(y-x) \,\mathrm{d}y \right) \mathrm{d}x = q\psi(0), \tag{15}$$

$$\lim_{t\to 0} \int_{\mathbb{R}} \psi(x) \left( \frac{1}{t} \int_{\mathbb{R}} q(t, x, y)(y-x)^2 \,\mathrm{d}y \right) \mathrm{d}x = \int_{\mathbb{R}} \psi(x) \,\mathrm{d}x. \tag{16}$$

The equality (16) identifies the diffusion coefficient, which is equal to $1$ here, while (15) allows us to identify the drift term of $X$ with $q\delta_0$. This is a generalization of the way of characterizing the drift term and the diffusion coefficient of a diffusion process, as presented for example in the book [71].

### 2.4. An explicit construction of the semi-group

We have constructed a semi-group $(Q_t)_{t>0}$ such that $Q_t\varphi(x)$ is a solution to the PDE (12). It is possible to look for a solution to (12) using an explicit computation.

Let $\varphi$ be a continuous, bounded function. Let us remark that

$$\nabla P_t \varphi(0) = \int_{\mathbb{R}} \frac{\sqrt{2}y}{\sqrt{\pi}t^{3/2}} \exp\left(-\frac{y^2}{2t}\right) \varphi(y)\,\mathrm{d}y$$

$$= \int_0^{+\infty} \frac{\sqrt{2}y}{\sqrt{\pi}t^{3/2}} \exp\left(-\frac{y^2}{2t}\right) (\varphi(y) - \varphi(-y))\,\mathrm{d}y$$

$$= \nabla P_t \widehat{\varphi}(0),$$

where $\widehat{\varphi}(x) = \mathbf{1}_{\{x \geq 0\}}(\varphi(x) - \varphi(-x))$. Let us also remark that $P_t \widehat{\varphi}(x) \xrightarrow[t \to 0]{} 0$ when $x \leq 0$. We set $\psi(x) = \widehat{\varphi}(-|x|)$. Thus, $P_t \psi(x)$ is solution to the heat equation on $\mathbb{R}_+^* \times \mathbb{R}^*$ with the initial condition equal to 0. Moreover, $\nabla P_t \psi(0+) = -\nabla P_t \varphi(0)$ and $\nabla P_t \psi(0-) = \nabla P_t \varphi(0)$. Thus, we are looking for a solution $T_t \varphi(x)$ written under the form

$$T_t \varphi(x) = P_t \varphi(x) + \lambda P_t \psi(x).$$

where $\lambda$ is chosen in order that $\alpha \nabla T_t \varphi(0+) = (1-\alpha)\nabla T_t \varphi(0-)$. Hence, it is always true if

$$\lambda = 2\alpha - 1 = q.$$

This result is a direct consequence of Proposition 1.

**Proposition 3.** *The semi-groups $(T_t)_{t>0}$ and $(Q_t)_{t>0}$ are equal.*

Hence, one gets the following expression for $Q_t$, which is more tractable than (14):

$$Q_t \varphi(x) = \int_{\mathbb{R}} p(t, x, y)\varphi(y)\,\mathrm{d}y$$

$$+ q \int_0^{+\infty} \frac{1}{\sqrt{2\pi t}} \exp\left(-\frac{(y + |x|)^2}{2t}\right) (\varphi(y) - \varphi(-y))\,\mathrm{d}y. \quad (17)$$

## 3. Using Dirichlet forms

The idea is to get rid of the singular first-order term by transforming $L$ into a symmetric operator in a proper Hilbert space.

### 3.1. *Weak solutions of PDE*

Here, we assume that $\alpha = (q+1)/2 \notin \{0, 1\}$. The problem with the parabolic PDE (12) is to deal with the transmission condition at 0 given by

$$\alpha \nabla u(t, 0+) = (1-\alpha)\nabla u(t, 0-). \quad (18)$$

So, an alternative way to consider (12) is to look for a weak solution of the following PDE

$$\begin{cases} u(t,x) \in \mathcal{C}(0,T;\mathrm{L}^2(\mathbb{R})) \cap \mathrm{L}^2(0,T;\mathrm{H}^1(\mathbb{R})), \\ \dfrac{\partial u(t,x)}{\partial t} = Au(t,x), \\ u(0,x) = \varphi(x) \in \mathrm{L}^2(\mathbb{R}), \end{cases} \tag{19}$$

where $A$ is the divergence form operator

$$A = \frac{1}{2a(x)} \nabla(a(x)\nabla\cdot) \text{ with } a(x) = \begin{cases} \alpha \text{ if } x \geq 0, \\ 1 - \alpha \text{ otherwise.} \end{cases} \tag{20}$$

By a weak solution to (19), we mean a function $u$ such that for all $\psi \in \mathcal{C}^\infty([0,T];\mathbb{R})$ with $\psi(T,x) = 0$, then, integrating formally (19) with respect to $\psi(t,x)\rho(x)\,\mathrm{d}x$ with $\rho(x) = 1/2a(x)$ and using integrations by parts,

$$\int_0^T \int_{\mathbb{R}} \frac{\partial \psi(t,x)}{\partial t} u(t,x)\rho(x)\,\mathrm{d}x\,\mathrm{d}t$$
$$+ \int_0^T \int_{\mathbb{R}} a(x)\frac{\partial u(t,x)}{\partial x}\frac{\partial \psi(t,x)}{\partial x}\,\mathrm{d}x\,\mathrm{d}t = -\int_{\mathbb{R}} \varphi(x)\psi(0,x)\rho(x)\,\mathrm{d}x.$$

As $u$ is smooth on the domain where the coefficient $a$ is smooth, it is easily deduced that $u$ is smooth on $(0,T]\times\mathbb{R}_+^*$ and $(0,T]\times\mathbb{R}_-^*$. Let us choose $\psi(t,x) = \psi_1(x)\psi_2(t)$ for two functions $\psi_1$ and $\psi_2$ smooth enough and such that $\psi_2(T) = 0$. Using an integration by parts on $\mathbb{R}_+^*$ and $\mathbb{R}_-^*$ and the freedom of choice of $\psi_1$ and $\psi_2$, we are led to (18). Moreover, one knows that $u(t,x)$ is continuous on $(0,T]\times\mathbb{R}$. Thus, the weak solution of (19) is also a solution of (12) and the converse is also true. The article [49] and the book [48, § III.13, p. 224] contain accounts on the properties of the solution of the transmission problem.

The domain $\mathrm{Dom}(A)$ of $A$ is the set of functions $f$ of $\mathrm{H}^1(\mathbb{R})$ such that $Af$ belongs to $\mathrm{L}^2(\mathbb{R})$. With this domain, $A$ is a self-adjoint (hence closed) operator with respect to the scalar product $\langle f,g\rangle_{\mathrm{L}^2(\mathbb{R};a)} = \int_{\mathbb{R}} f(x)g(x)a(x)\,\mathrm{d}x$ of $\mathrm{L}^2(\mathbb{R};a)$. It is also well known that a divergence-form operator as $A$ is the infinitesimal generator of a Feller semi-group $(\widehat{P}_t)_{t>0}$ with a density transition function $(t,x,y) \mapsto \Gamma(t,x,y)$ with respect to $a(x)\,\mathrm{d}x$, that is $\widehat{P}_t f(x) = \int_{\mathbb{R}} \Gamma(t,x,y)f(y)a(y)\,\mathrm{d}y$ for any continuous and bounded function $f$ on $\mathbb{R}$. Of course, $u(t,x) = \widehat{P}_t\varphi(x)$, and from the self-adjointness of $A$, $\Gamma(t,x,y) = \Gamma(t,y,x)$ for all $(t,x,y) \in \mathbb{R}_+^*\times\mathbb{R}^2$.

The Proposition 4 is a direct consequence of Proposition 1.

**Proposition 4.** *The semi-groups $(\widehat{P}_t)_{t>0}$ and $(Q_t)_{t>0}$ are equal.*

The next proposition follow from general results about divergence form operators: See [54, 78] for example.

**Proposition 5.** (i) *As $(\widehat{P})_{t>0}$ is a Feller semi-group, it is the generator of a strong Markov stochastic process which is continuous and conservative (which is then a SBM($\alpha$) thanks to Proposition 4).*

(ii) *There exists a constant $\delta \in (0,1)$ depending only on $\alpha$ such that the function $(t,x,y) \mapsto \Gamma(t,x,y)$ is $\delta/2$-Hölder continuous with respect to $t$ and $\delta$-Hölder continuous with respect to $(x,y)$.*

(iii) *There exist some constants $C_1, C_2, C_3, C_4 > 0$ depending only on $\alpha$ such that for all $(t,x,y) \in \mathbb{R}_+^* \times \mathbb{R}^2$,*

$$C_1 g(tC_2, x, y) \leq \Gamma(t, x, y) \leq C_3 g(tC_4, x, y), \tag{21}$$

*where $g(t,x,y) = \frac{1}{\sqrt{2\pi t}} \exp(-(x-y)^2/2t)$ is the heat kernel.*

The inequality (21) is called the *Aronson inequality*. It allows us to deduce some general results about the tightness of a family of processes generated by divergence form operators, or a bound on the probability to leave a sphere centered on $x$ prior to some time $T$ (see for example [54, 76, 78] for applications).

Of course, the transmission condition (18) implies that

$$\alpha \nabla_x \Gamma(t, 0-, y) = (1-\alpha) \nabla_x \Gamma(t, 0+, y), \ \forall (t,y) \in \mathbb{R}_+^* \times \mathbb{R}.$$

**A note on the adjoint of $A$ with respect to $L^2(\mathbb{R})$.** When one uses the scalar product $\langle f, g \rangle_{L^2(\mathbb{R})} = \int_{\mathbb{R}} f(x)g(x)\,\mathrm{d}x$, it is easily seen that the adjoint $A^*$ of $A$ is given by $A^* = \frac{1}{a}A(\frac{1}{a}\cdot)$. This implies that $\mathrm{Dom}(A^*) = \{ f \in L^2(\mathbb{R}) \mid f/a \in \mathrm{Dom}(A) \}$. Thus, any function $f$ in the domain of $A^*$ shall satisfies the condition

$$\alpha f(0-) = (1-\alpha) f(0+) \text{ and } \nabla f \text{ is continuous at } 0.$$

This also explains why the density $p(t,x,y)$ of the SBM($\alpha$) with respect to the Lebesgue measure is discontinuous as a function of $y$, in account to the Forward Kolmogorov equation. Of course, this can also be deduced from Proposition 5(ii).

### 3.2. Dirichlet forms

Let $(\mathcal{E}, \mathrm{Dom}(\mathcal{E}))$ be the quadratic form

$$\mathcal{E}(f, g) = \frac{1}{2} \int_{\mathbb{R}} a(x) \nabla f(x) \nabla g(x)\,\mathrm{d}x$$

for any $(f,g) \in \mathrm{Dom}(\mathcal{E}) = H^1(\mathbb{R})$. Let $\langle \cdot, \cdot \rangle_{L^2(a, \mathbb{R})}$ be the scalar product $\langle f, g \rangle_{L^2(a, \mathbb{R})} = \int_{\mathbb{R}} f(x)g(x)a(x)\,\mathrm{d}x$. Then $(\mathcal{E}, \mathrm{Dom}(\mathcal{E}))$ is the bilinear form associated to $(A, \mathrm{Dom}(A))$ by

$$\langle Af, g \rangle_{L^2(a, \mathbb{R})} = \mathcal{E}(f, g), \text{ for all } (f, g) \in \mathrm{Dom}(A) \times H^1(\mathbb{R}).$$

It is indeed easily checked that $(\mathcal{E}, \mathrm{Dom}(\mathcal{E}))$ is a *Dirichlet form* [59, 39], that is a symmetric bilinear form that is closed and contractive.

We can now give a third construction (however with somewhat complicated tools that we use here in a simple case) of a SBM($\alpha$).

**Proposition 6.** *The Dirichlet form $(\mathcal{E}, \mathrm{Dom}(\mathcal{E}))$ generates a continuous, Hunt (hence strong Markov) process, which is a SBM($\alpha$).*

*Proof.* This follows from standard results in the theory of Dirichlet forms since $(\mathcal{E}, \mathrm{Dom}(\mathcal{E}))$ is a regular, local Dirichlet form. Indeed, the infinitesimal generator of the process $(\mathcal{E}, \mathrm{Dom}(\mathcal{E}))$ generates is $(A, \mathrm{Dom}(A))$, hence it is a SBM($\alpha$). $\square$

*Remark* 1. In general, the process we construct using the previously invoked results on Dirichlet form is only defined for *quasi-every* starting point but not for any starting point. Yet in dimension one and under our assumptions on $\mathcal{E}$, any set of zero capacity is empty.

### 3.3. The Itô-Fukushima decomposition

The theory of Dirichlet form gives us another way to prove Theorem 1. We denote by $(X, \mathcal{F}_t, \mathbb{P}_x; t \geq 0, x \in \mathbb{R})$ the SBM($\alpha$) and by $(\theta_t)_{t \geq 0}$ its shift operator.

A *continuous additive functional* (CAF) is a continuous process $Y$ such that $Y_{t+s}(\omega) = Y_t(\theta_s \omega) + Y_s(\omega)$. A typical CAF is given by $\int_0^t f(X_s) \, \mathrm{d}s$ for some bounded, measurable function $f$.

A CAF $Y$ is characterized by a Radon measure $\nu$, called its *Revuz measure* by the following way: If $f$ is solution to

$$\langle f, \varphi \rangle_{\mathrm{L}^2(a, \mathbb{R})} + \mathcal{E}(f, \varphi) = \int_\mathbb{R} \varphi(x) \, \mathrm{d}\nu(x), \ \forall \varphi \in \mathcal{C}^1_\mathrm{c}(\mathbb{R}, \mathbb{R}),$$

then $\mathbb{E}_x\left[\int_0^{+\infty} e^{-t} \, \mathrm{d}Y_t\right]$ is a continuous[1] version of $f$.

Let $G_\alpha$ be the resolvent of $(\mathcal{E}, \mathrm{Dom}(\mathcal{E}))$, that is the operator giving the unique solution $u \in \mathrm{H}^1(\mathbb{R})$ to

$$\alpha \langle u, \varphi \rangle_{\mathrm{L}^2(a, \mathbb{R})} + \mathcal{E}(u, \varphi) = \langle f, \varphi \rangle_{\mathrm{L}^2(a, \mathbb{R})}, \ \forall \varphi \in \mathrm{H}^1(\mathbb{R}).$$

Using the resolvent of the Dirichlet form, we see immediately that if $\nu(\mathrm{d}x) = h(x) \, \mathrm{d}x$, the corresponding CAF is $\int_0^\cdot h(X_s)/a(X_s) \, \mathrm{d}s$.

We say that a CAF $Y$ is of *zero quadratic variation* if

$$e(Y) \stackrel{\mathrm{def.}}{=} \sup_{t \to 0} \frac{1}{2t} \int_\mathbb{R} a(x) \mathbb{E}_x\left[(Y_t - Y_0)^2\right] \mathrm{d}x = 0.$$

Of course, a CAF is said to be *locally of zero quadratic variation* if there exists a sequence $\{\tau^n\}_{n \in \mathbb{N}}$ of stopping times (with respect to the minimal filtration generated by $X$ satisfying the usual hypotheses) converging almost surely to $+\infty$ and such that $N_{\cdot \wedge \tau^n}$ is of zero quadratic variation.

**Theorem 2** (Itô-Fukushima decomposition). *Let $f$ be a function locally in $\mathrm{H}^1(\mathbb{R})$. Then there exists a local, square-integrable martingale $M^f$ with $\langle M^f \rangle_t = \int_0^f |\nabla f(X_s)|^2 \, \mathrm{d}s$, and a CAF $N^f$ locally of zero quadratic variation such that*

$$f(X_t) = f(X_0) + M_t^f + N_t^f.$$

---

[1] In the general case, it is only a quasi-continuous version, since a CAF may be only defined for any starting point except for those in a set of zero capacity. However, as noted in Remark 1, in dimension one and under our assumptions on $\mathcal{E}$, there are no such sets.

*Proof.* This theorem follows from the results in [39]. □

**Corollary 1.** *If $f$ belongs locally to $\mathrm{H}^2(\mathbb{R})$ (we choose a version of $f$ such that $f'$ is continuous), then*

$$N_t^f = \int_0^t f''(X_s)\,\mathrm{d}s + \beta_f \ell_t \ \text{ with } \beta_f = \frac{\alpha f'(0+) - (1-\alpha)f'(0-)}{2},$$

*where $(\ell_t)_{t\geq 0}$ is the CAF associated to the Revuz measure $\delta_0$.*

*If $f(x) = x$, then $N_t^f = \frac{(2\alpha-1)}{2}\ell_t$.*

*Remark 2.* Of course, up to some multiplicative constant, the CAF $\ell$ is the local time at $0$ of the process $X$ (see Proposition VI.45.10 in [74]).

*Proof.* Using a localization argument, we may assume that $f$ belongs to $\mathrm{H}^2(\mathbb{R})$, and we choose a version of $f$ such that $f'$ is continuous. By an integration by part, for any smooth function $\varphi$ with compact support,

$$\int_{\mathbb{R}} a(x)f'(x)\varphi'(x)\,\mathrm{d}x = -\int_0^{+\infty} \alpha f''(x)\varphi(x)\,\mathrm{d}x - \int_{-\infty}^0 (1-\alpha)f''(x)\varphi(x)\,\mathrm{d}x$$
$$- \alpha f'(0+)\varphi(0+) + (1-\alpha)f'(0-)\varphi(0-).$$

Hence,

$$\langle f, \varphi \rangle_{\mathrm{L}^2(a,\mathbb{R})} + \mathcal{E}(f,\varphi) = \int_{\mathbb{R}} \varphi(x)\nu(\,\mathrm{d}x)$$

with $\nu(\,\mathrm{d}x) = a(x)(f(x) - f''(x)/2)\,\mathrm{d}x - \beta_f \delta_0$. It follows that $f = G_1\nu$ and that

$$N_t^f = \int_0^t f(X_s)\,\mathrm{d}s - Y_t,$$

where $Y$ is a CAF associated to $\nu$. Yet the CAF associated to $\nu$ is $Y_t = \int_0^t (f(X_s) - f''(X_s)/2)\,\mathrm{d}s - \beta_f \ell_t$, which yields the result. □

### *3.4. Approximations of the coefficients*

We give now two results about the approximation of the coefficients of a divergence form operator.

Let $(\rho, a, b)$ be measurable functions from $\mathbb{R}$ to $\mathbb{R}$. We assume that there exist some constants $\lambda$ and $\Lambda$ for which

$$\lambda \leq a(x) \leq \Lambda, \ \lambda \leq \rho(x) \leq \Lambda \text{ and } |b(x)| \leq \Lambda. \tag{22}$$

Let $(L, \mathrm{Dom}(L))$ be the infinitesimal generator

$$L = \frac{1}{2\rho}\frac{\mathrm{d}}{\mathrm{d}x}\left(a\frac{\mathrm{d}}{\mathrm{d}x}\right) + b\frac{\mathrm{d}}{\mathrm{d}x}, \ \mathrm{Dom}(L) = \left\{ f \in \mathrm{H}^1(\mathbb{R}) \,\middle|\, Lf \in \mathrm{L}^2(\mathbb{R}) \right\}.$$

Under (22), $(L, \mathrm{Dom}(L))$ is the infinitesimal generator of a continuous, conservative strong Markov process (See [78, 54], ... for example).

**Proposition 7.** *Let $(\rho^n, a^n, b_n)_{n \in \mathbb{N}}$ be a family of measurable coefficients such that (22) holds uniformly and*

$$\frac{1}{a^n} \xrightarrow[n \to \infty]{\mathrm{L}^q(\mathbb{R})} \frac{1}{a(x)}, \quad \frac{1}{\rho^n} \xrightarrow[n \to \infty]{\mathrm{L}^q(\mathbb{R})} \frac{1}{\rho(x)} \ and \ \frac{b^n}{\rho^n a^n} \xrightarrow[n \to \infty]{\mathrm{L}^q(\mathbb{R})} \frac{b}{\rho a} \qquad (23)$$

*for some $1 \le q \le \infty$. Then the stochastic process $X^n$ generated by the operator $\frac{1}{2\rho^n} \frac{\mathrm{d}}{\mathrm{d}x} \left( a^n \frac{\mathrm{d}}{\mathrm{d}x} \right) + b^n \frac{\mathrm{d}}{\mathrm{d}x}$ converges in distribution to the process $X$ generated by $L$.*

*Proof.* Under (23), the convergence of the resolvent and the semi-group of $L^n$ to that of $L$ follows from results in [99, 100] for example. Hypothesis (22) implies that one may compare the density transition functions of the processes $X^n$ and $X$ uniformly with respect to the Gaussian kernel, and that these functions are locally Hölder continuous in space and time. Thus, this is sufficient to deduce the weak convergence of $X^n$ to $X$: see [54, 76] for example. □

**Corollary 2.** *Let $a^n$, $\rho^n$ and $b^n$ be such that $(\rho^n, a^n, b^n)_{n \in \mathbb{N}}$ satisfies (22) and*

$$(\rho^n(x), a^n(x), b^n(x)) \xrightarrow[n \to \infty]{} (\rho(x), a(x), b(x)) \ almost \ everywhere,$$

*then the conclusion of Proposition 7 is true.*

In order to approximate the SBM by some classical SDEs, it is then possible to use Proposition 7 or Corollary 2 with the operator $A$ defined in (20). We will see in the next section an application of this result to identify to different constructions of the SBM. We give now another application, in which the SBM is constructed as a renormalized limit of some SDE with drift. The following result has first appeared in [75] (See also Theorem 14 and Remark 15 below for a restatement of this result in the more general context of SDEs with local time).

**Proposition 8.** *Let $X$ be the solution to the SDE*

$$X_t = B_t + \int_0^t b(X_s) \, \mathrm{d}s \qquad (24)$$

*where $B$ is a Brownian motion and $b \in \mathrm{L}^1(\mathbb{R}; \mathbb{R})$. If $X_t^n = n^{-1} X_{tn^2}$, then $X^n$ converges in distribution to a SBM($\alpha$) with $\alpha = e^\kappa/(1 + e^\kappa)$ and $\kappa = 2 \int_{-\infty}^{+\infty} b(y) \, \mathrm{d}y$.*

*Proof.* In dimension one, the infinitesimal generator $A = \frac{1}{2\rho} \frac{\mathrm{d}}{\mathrm{d}x} \left( a \frac{\mathrm{d}}{\mathrm{d}x} \right) + b \frac{\mathrm{d}}{\mathrm{d}x}$ may be written

$$A = \frac{e^{-\Phi}}{2\rho} \frac{\mathrm{d}}{\mathrm{d}x} \left( a e^\Phi \frac{\mathrm{d}}{\mathrm{d}x} \right) \ \mathrm{with} \ \Phi(x) = \int_\zeta^x \frac{2b(y)}{a(y)\rho(y)} \, \mathrm{d}y$$

for any $\zeta \in \mathbb{R}$. As $b$ belongs to $\mathrm{L}^1(\mathbb{R}; \mathbb{R})$, then one can set $\zeta = -\infty$. Assume that $a = \rho = 1$. The process $(X, \mathbb{P}_x)$ generated by $A$ is then the solution to (24). The infinitesimal generator of $X^n$ is $\frac{1}{2} \triangle + n b(\cdot n) \frac{\mathrm{d}}{\mathrm{d}x}$.

It follows that, if $\Phi^n(x) = 2 \int_{-\infty}^x nb(ny) \, \mathrm{d}y$, then $\Phi^n(x) = 2 \int_{-\infty}^{nx} b(y) \, \mathrm{d}y$ and converges to $\kappa$ when $x > 0$ and to $0$ when $x < 0$. Thus with (20) where $a(x)$ has been divided by $(1-\alpha)$ and Corollary 2, it is easily checked that $X^n$ converges in distribution to the Skew Brownian motion of parameter $\alpha = e^\kappa/(1+e^\kappa)$. $\square$

## 4. Using the scale function and the speed measure

### 4.1. The scale function and the speed measure

Studying a Markov process $(X, \mathbb{P}_x)_{x \in \mathbb{R}}$ with infinitesimal generator $(L, \mathrm{Dom}(L))$ is much simpler in dimension one, since its behaviour can generally be described by two unique (up to multiplicative and additive constants) strictly increasing functions $S$ and $V$. On this description, see the books [12, 44, 74],...

The function $S$, defined on the state space of the process, is called the *scale function* and satisfies

$$\mathbb{P}_y[\, \tau_x < \tau_z \,] = \frac{S(z) - S(y)}{S(x) - S(z)} \text{ for all } x < y < z,$$

where $\tau_x = \inf \{\, t > 0 \,|\, X_t = x \,\}$ is the first hitting time of the point $x$ in the state space. The function $S$ is harmonic for $L$, and the scale function of the Brownian motion is $S(x) = x$. A process for which $S(x) = x$ is said to be *on its natural scale*. Of course, if $S^{-1}$ is the right-continuous inverse of $S$, $\widehat{X} = S(X)$ is also a Markov process and is on its natural scale.

The function $V$ is a primitive of a measure $m$ on the state space of the process, which is called the *speed measure*. Indeed, this measure characterizes the average exit time from any given interval for the process. If $a < x < b$, then

$$\mathbb{E}_x[\, \tau_{(a,b)} \,] = \int_a^b G_{(a,b)}(x,y) m(\mathrm{d}y)$$

where $G_{(a,b)}$ is the Green function

$$G_{(a,b)}(x,y) = \begin{cases} \frac{2(S(x) - S(a))(S(b) - S(y))}{S(b) - S(a)} & \text{if } x \leq y, \\ \frac{2(S(y) - S(a))(S(b) - S(x))}{S(b) - S(a)} & \text{if } x > y \end{cases}$$

and $\tau_{(a,b)} = \inf \{\, t > 0 \,|\, X_t \notin (a,b) \,\}$. Any Markov process on its natural scale may be written as $B_{T(\cdot)}$, where $B$ is a Brownian motion and $T(\cdot)$ is a random time change computed from $B$ and $m$ (or equivalently $V$). If $m$ is absolutely continuous with respect to the Lebesgue measure, then the time spent by the process at any point has zero Lebesgue measure.

*Remark* 3. The process $\widehat{X}$ is the scale function $\widehat{S}(x) = x$ and the integrated speed measure $\widehat{V}(x) = V \circ S^{-1}$. If the speed measure $m$ of $X$ has a density $\varphi(x)$ with respect to the Lebesgue measure, then the speed measure $\widehat{m}$ has a density $\varphi(S^{-1}(x))/S'(S^{-1}(x))$ with respect to the Lebesgue measure.

The boundary conditions are also coded by $S$ and $m$, but from now, we only deal with processes whose state-space is $\mathbb{R}$ and that does not explode.

The infinitesimal generator $(L, \mathrm{Dom}(L))$ can be described by

$$L = \frac{1}{2} \frac{\mathrm{d}}{\mathrm{d}V} \frac{\mathrm{d}}{\mathrm{d}S}, \tag{25}$$

where

$$\frac{\mathrm{d}f}{\mathrm{d}\varphi}(x) = \lim_{\varepsilon \to 0} \frac{f(x + \varepsilon) - f(x)}{\varphi(x + \varepsilon) - \varphi(x)}$$

for $\varphi = S$ or $V$. At a point $x$ where $\varphi$ is differentiable, $\frac{\mathrm{d}f}{\mathrm{d}\varphi}(x) = f'(x)/\varphi'(x)$.

The following theorem provides a way to prove the convergence of stochastic processes by looking at their scale functions and speed measures. It allows us to identify the process constructed by its speed measure and its scale function with a process constructed by Dirichlet forms or PDEs results.

**Theorem 3** ([36]). *Let $(X, \mathbb{P}_x; x \in \mathbb{R})$ be a conservative, continuous process defined by the functions $(S, V)$, and let $(X, \mathbb{P}_x^n; x \in \mathbb{R})_{n \in \mathbb{N}}$ be a family of conservative, continuous processes defined by $(S^n, V^n)$. If for all $x \in \mathbb{R}$, $S^n(x) \xrightarrow[n \to \infty]{} S(x)$ and $V^n(x) \xrightarrow[n \to \infty]{} V(x)$, then $\mathbb{P}_x^n$ converges weakly to $\mathbb{P}_x$ in the space of continuous functions for any starting point $x$.*

### *4.2. The Skew Brownian motion*

Comparing (25) and (20) suggests to construct the process $(X, \mathbb{P}_x)_{x \in \mathbb{R}}$ with scale function and integrated speed measures

$$S(x) = \begin{cases} \alpha^{-1}x & \text{if } x \geq 0, \\ (1-\alpha)^{-1}x & \text{if } x < 0, \end{cases} \text{ and } V(x) = \begin{cases} \alpha x & \text{if } x \geq 0, \\ (1-\alpha)x & \text{if } x < 0. \end{cases} \tag{26}$$

*Remark* 4. If $(S, V)$ are the scale function and the integrated speed measure of a diffusion process, then $(\kappa S + \lambda, \kappa^{-1}V + \lambda)$ for $\kappa > 0$ and $\lambda \in \mathbb{R}$ are also possible scale functions and integrated speed measure.

The question is to know if $(X, \mathbb{P}_x)_{x \in \mathbb{R}}$ is equal to the process previously constructed with Feller semi-groups.

**Proposition 9.** *The scale function $S$ and the integrated speed measure $V$ of the SBM($\alpha$) are given by (26).*

*Proof.* Let $(a^n, \rho^n, b^n)_{n \in \mathbb{N}}$ be a family of measurable functions satisfying (22). Let us consider

$$A^n = \frac{1}{2\rho^n} \frac{\mathrm{d}}{\mathrm{d}x} \left( a^n \frac{\mathrm{d}}{\mathrm{d}x} \right) + b^n \frac{\mathrm{d}}{\mathrm{d}x}. \tag{27}$$

If $a^n$, $\rho^n$ and $b^n$ are smooth enough, then it is known that $A^n$ is the infinitesimal generator of a continuous stochastic process $(X, \mathbb{P}_x^n)_{x \in \mathbb{R}}$ that can be described

either by the functions

$$S^n(x) = \int_0^x \frac{\exp(-h^n(y))}{a^n(y)}\,\mathrm{d}y \text{ and } V^n(x) = \int_0^x \frac{\exp(-h^n(y))}{\rho^n(y)}\,\mathrm{d}y$$

$$\text{with } h^n(x) = 2\int_0^x \frac{b^n(y)}{\rho^n(y)a^n(y)}\,\mathrm{d}y,$$

or by its semi-group, which is a Feller semi-group. Hence, if almost everywhere, $a^n(x) \xrightarrow[n\to\infty]{} a(x)$, $\rho^n(x) \xrightarrow[n\to\infty]{} 1/a(x)$ and $b^n(x) \xrightarrow[n\to\infty]{} 0$, then Proposition 7 together with Theorem 3 prove that the process generated by $A$ *via* its Feller semi-group and the one by its scale function and speed measure are equal in distribution. □

## 5. The SBM as solution of a SDE with local time

Combining the description of the Skew Brownian motion together with the Itô-Tanaka formula allows us to strengthen the semi-martingale decomposition in Theorems 1 and 2 and to express the Skew Brownian motion as the strong solution of some SDE with its local time.

### 5.1. The Itô-Tanaka formula

Let $f$ is the function from $\mathbb{R}$ to $\mathbb{R}$ which is the difference of two convex functions. Then $f'_r(x) = \lim_{\varepsilon\to 0,\varepsilon>0}\epsilon^{-1}(f(x+\varepsilon)-f(x))$ and $f'_\ell(x) = \lim_{\varepsilon\to 0,\epsilon<0}\epsilon^{-1}(f(x+\varepsilon)-f(x))$ exists for almost every $x$. In addition, there exists a signed measure $\mu$, called the *second derivative measure*, such that

$$\int_{\mathbb{R}} g(x)\mu(\,\mathrm{d}x) = -\int_{\mathbb{R}} g'(x)f'_\ell(x)\,\mathrm{d}x. \tag{28}$$

for any piecewise $\mathcal{C}^1$ function with compact support on $\mathbb{R}$. If $f$ has a second derivative, then $f''$ is the density of $\mu$ with respect to the Lebesgue measure.

Let $X$ be a real-valued semi-martingale. Then there exists a process $(L_t^{x-}(X))_{t\geq 0,x\in\mathbb{R}}$, called the *left local time*, such that for any function $f$ as above,

$$f(X_t) = f(x) + \int_0^t f'_\ell(X_s)\,\mathrm{d}X_s + \frac{1}{2}\int_{\mathbb{R}} L_t^{x-}(X)\mu(\,\mathrm{d}x). \tag{29}$$

This is the *Itô-Tanaka formula*. The left local time is continuous, has finite variation and satisfies $L_t^{x-}(X) = \int_0^t \mathbf{1}_{\{X_s=x\}}\,\mathrm{d}L_s^{x-}$. Indeed for $t\geq 0$ and $x\in\mathbb{R}$, the left local time $L_t^{x-}(X)$ may be defined by

$$\frac{1}{2}L_t^{x-}(X) = \int_0^t \mathbf{1}_{\{X_s>x\}}\,\mathrm{d}X_s - (X_t-x)^+. \tag{30}$$

*Remark* 5. Here, we use a normalisation of the local time which is different from the one in [46].

We also set for $x \in \mathbb{R}$,

$$\frac{1}{2}L_t^{x+}(X) = \int_0^t \mathbf{1}_{\{X_s < x\}} \, \mathrm{d}X_s - (X_t - x)^-, \tag{31}$$

which we call the *right local time*. For $Y_t = -X_t$, we have that

$$(Y_t + x)^+ = -\int_0^t \mathbf{1}_{\{X_s < -x\}} \, \mathrm{d}X_s + \frac{1}{2}L_t^{(-x)-}(Y)$$

$$= -(X_t - x)^- = -\int_0^t \mathbf{1}_{\{X_s < x\}} \, \mathrm{d}X_s - \frac{1}{2}L_t^{x+}(X).$$

It follows that $L_t^{x+}(X) = -L_t^{(-x)-}(Y)$. Now, if $g(x) = f(-x)$, then $g_r'(x) = -f_\ell'(-x)$ and $g_\ell'(x) = -f_r'(-x)$. As $f_r'(x)$ differs from $f_\ell'(x)$ on a countable set [46, Problem 6.21, p. 213], the second derivative measure of $g$ is equal to $\mu(-\,\mathrm{d}x)$. Hence, with the Itô-Tanaka formula,

$$g(Y_t) = g(-x) - \int_0^t f_r'(-Y_s) \, \mathrm{d}Y_s - \frac{1}{2}\int_{\mathbb{R}} \mu(-\,\mathrm{d}x)L_t^{x-}(Y)$$

$$= f(x) + \int_0^t f_r'(X_s) \, \mathrm{d}X_s + \frac{1}{2}\int_{\mathbb{R}} \mu(\,\mathrm{d}x)L_t^{x+}(X).$$

As $g(Y_t) = f(X_t)$, summing this expression with (29), we get

$$f(X_t) = f(x) + \int_0^t \frac{1}{2}(f_r'(X_s) + f_\ell'(X_s)) \, \mathrm{d}X_s + \frac{1}{2}\int_0^t \mu(\,\mathrm{d}x)L_t^x(X), \tag{32}$$

where

$$L_t^0(x) = \frac{L_t^{0-}(X) + L_t^{0+}(X)}{2} \tag{33}$$

is the *symmetric local time*. Formula (32) is the *symmetric Itô-Tanaka formula*. Let us note that if $M$ is the martingale part of $X$, then $\int_0^t \mathbf{1}_{\{X_s = x\}} \, \mathrm{d}M_s$ is equal to 0, as it follows easily from computing the expectation of the brackets of the integral.

*Remark* 6. We have defined the left, right and symmetric local time of a diffusion process $X$ at 0 in view of using the Itô-Tanaka formula. An alternative construction is

$$L_t^0(X) = \lim_{\varepsilon \to 0} \frac{1}{2\varepsilon}\int_0^t \mathbf{1}_{\{X_s \in [-\varepsilon, \varepsilon]\}} \, \mathrm{d}s$$

for the symmetric local time at 0, and

$$L_t^{0-}(x) = \lim_{\varepsilon \to 0} \frac{1}{\varepsilon}\int_0^t \mathbf{1}_{\{X(s) \in [-\varepsilon, 0]\}} \, \mathrm{d}s \text{ and } L_t^{0+}(x) = \lim_{\varepsilon \to 0} \frac{1}{\varepsilon}\int_0^t \mathbf{1}_{\{X(s) \in [0, \varepsilon]\}} \, \mathrm{d}s$$

for the left and right local time at 0. This is a consequence of the *occupation time formula* [46, Theorem 7.1(iii)] which asserts that if $f$ is a measurable, bounded function, then

$$\int_{\mathbb{R}} f(x)L_t^x(X) \, \mathrm{d}x = \int_0^t f(X_s) \, \mathrm{d}\langle X \rangle_s \tag{34}$$

for all $t \geq 0$. This formula follows easily from the Itô formula when $f$ is smooth enough, and then by a density argument.

### 5.2. The SDE the SBM solves

Using the scale function and the speed measure, it is then easy to construct the SBM as the solution of some SDE. For this, let us consider first the SDE

$$Y_t = y + \int_0^t \sigma(Y_s)\, \mathrm{d}B_s \text{ with } \sigma(y) = \frac{1}{\alpha}\mathbf{1}_{\{y \geq 0\}} + \frac{1}{1-\alpha}\mathbf{1}_{\{y < 0\}}, \quad (35)$$

where $B$ is a one-dimensional Brownian motion and $y \in \mathbb{R}$. According to a result due to S. Nakao [63] (see also [51, 52]), this SDE has a unique strong solution.

Thus, let $B$ be a standard Brownian motion with its natural filtration (transformed to satisfy the standard hypotheses). The strong solution to (35) with respect to $B$ is also given by $Y_t = y + B_{T(t)}$ where $T(t)$ is the random time change defined by

$$t = \int_0^{T(t)} \frac{\mathrm{d}s}{\sigma(B_s)^2}. \quad (36)$$

The inverse $r(x)$ of $S(x)$ is $r(x) = \alpha x$ if $x \geq 0$ and $(1-\alpha)x$ if $x < 0$. With the symmetric Itô-Tanaka formula (32),

$$r(Y_t) = r(y) + \int_0^t \alpha \mathbf{1}_{\{Y_s > 0\}}\, \mathrm{d}Y_s + \int_0^t (1-\alpha)\mathbf{1}_{\{Y_s < 0\}}\, \mathrm{d}Y_s + \frac{2\alpha-1}{2}L_t^0(Y),$$

since $\int_0^t \mathbf{1}_{\{Y_s=0\}}\, \mathrm{d}Y_s = 0$. Now, let us remark that

$$\int_0^t \alpha \mathbf{1}_{\{Y_s > 0\}}\, \mathrm{d}Y_s + \int_0^t (1-\alpha)\mathbf{1}_{\{Y_s < 0\}}\, \mathrm{d}Y_s = \int_0^t \mathbf{1}_{\{B_s \neq 0\}}\, \mathrm{d}B_s = B_t$$

since $\int_0^t \mathbf{1}_{\{B_s=0\}}\, \mathrm{d}B_s = 0$. Thus, if $X_t = r(Y_t) = S^{-1}(B_{T(t)})$, this proves that $S$ is the scale function of $X$ and $V$ is the integrated speed measure of $X$ (The function $V$ is identified through (36) and Remark 3 with the help of Theorem 16.84 in [12] for example).

Hence,

$$X_t = B_t + \frac{2\alpha-1}{2}L_t^0(Y)$$

for any $t \geq 0$. It remains to express $L_t^0(Y)$ in terms of $L_t^0(X)$. Again with the symmetric Itô-Tanaka formula,

$$|X_t| = \int_0^t \mathrm{sgn}(X_s)\, \mathrm{d}X_s + L_t^0(X) = \int_0^t \mathrm{sgn}(X_s)\, \mathrm{d}B_s + L_t^0(X),$$

where $\mathrm{sgn}(x) = 1$ if $x > 0$, $-1$ if $x < 0$ and $0$ if $x = 0$. On the other hand, the second derivative measure of $|r(x)|$ is the Dirac measure $\delta_0$ at $0$. Hence, since $\mathrm{sgn}(X_t) = \mathrm{sgn}(Y_t)$,

$$|r(Y_t)| = \int_0^t \mathrm{sgn}(Y_s) r'(Y_s)\, \mathrm{d}Y_s + \frac{1}{2}L_t^0(Y) = \int_0^t \mathrm{sgn}(X_s)\, \mathrm{d}B_s + \frac{1}{2}L_t^0(Y).$$

From the uniqueness of the decomposition of $X$ as a semi-martingale that $L_t^0(X) = \frac{1}{2}L_t^0(Y)$. Thus, $X$ is the strong solution to

$$X_t = x + B_t + \beta L_t^0(X) \text{ with } \beta = 2\alpha - 1, \tag{37}$$

where $L_t^0(X)$ is the symmetric local time of $X$ at 0 and $x = r(y)$. As $r$ is one-to-one and since the solution to (35) is unique, one gets also that $X$ is unique.

**Theorem 4** ([41]). *The SDE* (37) *has a unique strong solution if and only if* $\alpha \in [0, 1]$. *If a strong solution exists, then it is the Skew Brownian motion of parameter* $\alpha$.

The construction of the SDE (37) relies on the fact that $S$ is one-to-one, $S'$ is constant on $\mathbb{R}_+$ and $\mathbb{R}_-^*$ and $S(\mathbb{R}_\pm) \subset \mathbb{R}_\pm$. Then (35) is equivalent to (37). We then also get that weak existence and weak uniqueness for (37), as the SDE (35) has also a unique weak solution. This will be used below to deal with the martingale problem.

*Remark* 7. In [52], J.-F. Le Gall extended results on strong existence and uniqueness to SDEs of type $\mathrm{d}X_t = \sigma(X_t)\,\mathrm{d}B_t + \int_{\mathbb{R}} \nu(\,\mathrm{d}x)\,\mathrm{d}L_t^x(X)$ under rather general conditions on the coefficient $\sigma$ and the measure $\nu$. In this article, the reader will also find another way to prove that the diffusion process generated by (27) also converges to the solution of (37).

More recently, R. Bass and Z.-Q. Chen also considered this kind of equation in [9] and extended some of the results of [52].

Section 11.8 contains a short account of this theory.

### *5.3. About the left, right and symmetric local time*

As $t \mapsto L_t^0(X)$ is a continuous additive functional that increases only when $t$ belongs to $\mathcal{Z} = \overline{\{s \geq 0 \,|\, X_s = 0\}}$, any continuous additive functional $t \mapsto \eta_t$ with the same property is almost surely equal to $C_\eta L^0(X)$, where $C_\eta$ is a constant that depends only on $\eta$ (see Proposition VI.45.10 in [74]). Hence, $L^0(X)$, $L^{0-}(X)$ and $L^{0+}(X)$ are all proportional.

It follows from (30) and (31) that

$$L_t^{0+}(X) - L_t^{0-}(X) = 2 \int_0^t \mathbf{1}_{\{X_s = 0\}}\,\mathrm{d}X_s. \tag{38}$$

As the quadratic variation of $\int_0^{\cdot} \mathbf{1}_{\{X_s = 0\}}\,\mathrm{d}B_s$ is equal to 0, this term is also equal to 0. But $\int_0^t \mathbf{1}_{\{X_s = 0\}}\,\mathrm{d}L_s^0(X) = L_t^0(X)$. It follows that

$$L_t^{0+}(X) = 2\alpha L_t^0(X) \text{ and } L_t^{0-}(X) = 2(1-\alpha)L_t^0(X). \tag{39}$$

For all $t \geq 0$, $x \mapsto L_t^{x+}(X)$ is right-continuous with a left-limit at each point. Moreover, $L_t^{x-}(X) = \lim_{\varepsilon \to 0, \ \varepsilon < 0} L_t^{(x-\varepsilon)+}(X)$. Outside 0, the process $X$ behaves like the Brownian motion and thus $x \mapsto L_t^{x+}(X)$ is continuous on $\mathbb{R} \setminus \{0\}$. With equalities (39), this proves the following theorem, initially due to J. Walsh.

**Theorem 5** ([87]). *Unless for* $\alpha = 1/2$ *(the Brownian case), the map* $x \mapsto L_t^{x+}(X)$ *is discontinuous at 0 and continuous elsewhere.*

## 6. The martingale problem

Of course, one may wish to describe the distribution of the SBM($\alpha$) with the help of the martingale problem (see [46, 79] among many other books).

Due to the presence of the local time in the SDE describing the SBM($\alpha$), this is not the most suitable construction. Yet we will be able to define a martingale problem and show that it is well posed.

Let us denote by $D(\alpha)$ the set of continuous, bounded functions $f$ on $\mathbb{R}$ with two bounded derivatives $f'$ and $f''$ on $\mathbb{R}^*$ such that $f''(0+)$ and $f''(0-)$ exist, and $\alpha f'(0+) = (1-\alpha) f'(0-)$. Hence, $f$ is the difference of two convex functions and its second generalized derivative is

$$\mu(\,\mathrm{d}x) = f''(x)\,\mathrm{d}x + \frac{1-2\alpha}{\alpha} f'(0-)\delta_0.$$

The SBM($\alpha$) is the diffusion process $(X, \mathcal{F}_t, \mathbb{P}_x; t \geq 0, x \in \mathbb{R})$ which is solution to $\mathrm{d}X_t = \mathrm{d}B_t + \beta L_t^0(X)$ with $\beta = 2\alpha - 1$. Using the symmetric Itô-Tanaka formula (32),

$$
\begin{aligned}
f(X_t) &= f(x) + \int_0^t f'(X_s)\,\mathrm{d}B_s + \frac{1}{2}\int_{\mathbb{R}} \mu(\,\mathrm{d}x) L_t^x(X) \\
&\qquad + \frac{\beta}{2}(f'(0+) + f'(0-)) L_t^0(X) \\
&= f(x) + \int_0^t f'(X_s)\,\mathrm{d}B_s + \frac{\beta}{2\alpha} f'(0-) L_t^0(X) \\
&\qquad + \frac{1}{2}\int_{\mathbb{R}} f''(x) L_t^x(X)\,\mathrm{d}x + \frac{1-2\alpha}{2\alpha} f'(0-) L_t^0(X).
\end{aligned}
$$

With the occupation time formula (34), since $\langle X \rangle_t = t$, $\int_{\mathbb{R}} f''(x) L_s^{x-}(X)\,\mathrm{d}x = \int_0^t f''(X_s)\,\mathrm{d}s$. As $\beta = 2\alpha - 1$, we obtain that

$$\forall f \in D(\alpha),\ f(X_t) = f(x) + \int_0^t f'(X_s)\,\mathrm{d}B_s + \frac{1}{2}\int_0^t f''(X_s)\,\mathrm{d}s$$

for all $x \in \mathbb{R}$. Thus,

$$M_t^f = f(X_t) - f(x) - \frac{1}{2}\int_0^t f''(X_s)\,\mathrm{d}s \tag{40}$$

is a $(\mathcal{F}_t)_{t \geq 0}$-martingale under $\mathbb{P}_x$.

We denote by $(X_t)_{t \geq 0}$ the canonical process on $\mathcal{C}(\mathbb{R}_+; \mathbb{R})$. Let us define by $(\mathcal{B}_t)_{t \geq 0}$ the filtration $\mathcal{B}_t = \sigma(X_s; s \leq t)$.

Thus, let us now state the martingale problem.

**Definition 1.** A probability measure $\mathbb{Q}$ on $(\mathcal{C}(\mathbb{R}_+; \mathbb{R}), \mathrm{Bor}(\mathcal{C}(\mathbb{R}_+; \mathbb{R})))$ for which for all $0 \leq s \leq t$,

$$\mathbb{E}^{\mathbb{Q}}\left[ f(X_t) - f(X_s) - \frac{1}{2}\int_s^t f''(X_r)\,\mathrm{d}r \,\bigg|\, \mathcal{B}_s \right] = 0,\ \forall f \in D(\alpha),$$

is said to be a *solution of the martingale problem associated to $D(\alpha)$*.

**Proposition 10.** *The distribution $\mathbb{P}_x$ of the SBM($\alpha$) is the unique solution of the martingale problem associated to $D(\alpha)$ satisfying $\mathbb{P}_x[\,X_0 = x\,] = 1$ (In other words, the martingale problem is well posed).*

*Proof.* We have already seen that the SBM($\alpha$) is solution to the martingale problem associated to $D(\alpha)$, except for the replacement of the filtration of the Brownian motion by $(\mathcal{B}_t)_{t\geq0}$. Yet, this last point may be treated as usual (see for example Remark 5.4.16 in [46, p. 320]).

Let $\mathbb{Q}$ be a solution of the martingale problem with $\mathbb{Q}[\,X_0 = x\,] = 1$.

Let us consider the scale function $S$ of the SBM($\alpha$): $S(x) = \alpha^{-1}x$ if $x \geq 0$ and $S(x) = (1-\alpha)^{-1}x$ if $x < 0$. Then $S''(x) = 0$ for $x \neq 0$. The function $S$ is not bounded, but with a localization argument, one gets easily that

$$S(X_t) = S(x) + M_t, \quad \mathbb{Q}\text{-a.s.},$$

where $M$ is a local martingale. Let us compute its brackets: with the Itô formula applied to $x \mapsto x^2$,

$$S(X_t)^2 = S(x)^2 + 2\int_0^t S(X_s)\,\mathrm{d}M_s + \langle M \rangle_t.$$

Now, if $U(x) = S^2(x)$, then $U$ also satisfies $\alpha U'(0+) = (1-\alpha)U'(0-)$ since $U'(0+) = U'(0-) = 0$. In addition, $U''(x) = 2S'(x)^2$. Also with a localization argument,

$$U(X_t) = S(x)^2 + N_t + \int_0^t S'(X_s)^2\,\mathrm{d}s, \quad \mathbb{Q}\text{-a.s.},$$

where $N$ is a local martingale. Since

$$M_t^2 = (S(X_t) - S(x))^2 = 2\int_0^t S(X_s)\,\mathrm{d}M_s + \langle M \rangle_t - S(x)M_t$$

$$= \int_0^t S'(X_s)^2\,\mathrm{d}s - S(x)N_t,$$

one deduces easily that $\langle M \rangle_t = \int_0^t S'(X_s)^2\,\mathrm{d}s$. From the representation theorem (see for example Theorem 4.2 in [46, p. 170]), there exists an extension $(\widetilde{\Omega}, \widetilde{\mathcal{F}}, \widetilde{\mathbb{Q}})$ of the probability space $(\mathcal{C}(\mathbb{R}_+; \mathbb{R}), \mathrm{Bor}(\mathcal{C}(\mathbb{R}_+; \mathbb{R})), \mathbb{Q})$ as well as a filtration $(\widetilde{\mathcal{F}}_t)_{t\geq0}$ and a $(\widetilde{\mathbb{Q}}, (\widetilde{\mathcal{F}}_t)_{t\geq0})$-Brownian motion $\widetilde{B}$ on this space such that $M_t = \int_0^t \rho(s)\,\mathrm{d}\widetilde{B}_s$, where $\rho$ is $(\widetilde{\mathcal{F}}_t)_{t\geq0}$-adapted. Besides, $\widetilde{\mathbb{Q}}$-almost surely, $\rho^2 = S'(X.)^2$. Hence, setting $W_t = \int_0^t \mathrm{sgn}(\rho(s)S'(X_s)^2)\,\mathrm{d}\widetilde{B}_s$, one gets that, $\widetilde{\mathbb{Q}}$-almost surely, $M_t = \int_0^t S'(X_s)\,\mathrm{d}W_s$ for any $t \geq 0$.

Since $S'$ is constant on $\mathbb{R}_+$ and $\mathbb{R}_-^*$ and $S(\mathbb{R}_\pm) \subset \mathbb{R}_\pm$, one can set $Y_t = S(X_t)$ and then $Y_t = S(x) + \int_0^t S'(Y_s)\,\mathrm{d}W_s$ for any $t \geq 0$, $\widetilde{\mathbb{Q}}$-almost surely. Thus, $(Y, (\widetilde{\mathcal{F}}_t)_{t\geq0}, \widetilde{\mathbb{Q}})$ is a weak solution to (35), which is known to be unique. Thus $S^{-1}(Y) = X$ is the SBM($\alpha$). $\qquad\square$

## 7. Decomposition of the excursions' measure

Another construction is the following: Let $Y$ be a Reflected Brownian motion, and let $\mathcal{Z} = \overline{\{\, s \geq 0 \,|\, Y_s = 0 \,\}}$ be the closure of its zeros. The Lebesgue measure of $\mathcal{Z}$ is zero, but this set cannot be ordered. However, the set $\mathbb{R}_+ \setminus \mathcal{Z}$ can be decomposed as a countable union $\cup_{n \in \mathbb{N}} J_n$ of intervals $J_n$. Each interval $J_n$ corresponds to some excursion of $Y$, that is, if $J_n = (\ell_n, r_n)$,

$$Y_t > 0 \text{ for } t \in (\ell_n, r_n) \text{ and } Y_{\ell_n} = Y_{r_n} = 0.$$

Fix $\alpha \in [0, 1]$. At each $J_n$, we associate a Bernoulli random variable $e_n$ which is independent from any other random variables (and the Reflected Brownian motion $Y$) and such that $\mathbb{P}[\, e_n = 1 \,] = \alpha$ and $\mathbb{P}[\, e_n = -1 \,] = 1 - \alpha$. Let $X$ be the process given by

$$X_t = e_n Y_t \text{ if } t \in J_n.$$

**Theorem 6** ([44]). *The process $X$ is a Skew Brownian motion of parameter $\alpha$.*

*Proof.* Obviously, the process $X$ behaves like a Brownian motion on any time interval $J_n$. By Remark 4, its scale function $S$ and its integrated speed measure $V$ are of type

$$S(x) = \begin{cases} \gamma^{-1} x & \text{if } x \geq 0, \\ (1-\gamma)^{-1} x & \text{if } x \leq 0 \end{cases} \text{ and } V(x) = \begin{cases} \gamma x & \text{if } x \geq 0, \\ (1-\gamma) x & \text{if } x \leq 0 \end{cases}$$

where the constant $\gamma$ has to be specified. If $\tau_1$ (resp. $\tau_{-1}$) are the first time the process $X$ reaches $1$ (resp. $-1$),

$$\mathbb{P}_0[\, \tau_1 < \tau_{-1} \,] = \frac{S(-1) - S(0)}{S(-1) - S(1)} = \gamma.$$

On the other hand, using the decomposition by excursions and the independence of the $e_n$'s with respect to $Y$,

$$\mathbb{P}_0[\, \tau_1 < \tau_{-1} \,] = \sum_{n \in \mathbb{N}} \mathbb{P}_0[\, \tau_1 \wedge \tau_{-1} \in J_n; \; X_t \geq 0 \text{ on } J_n \,]$$

$$= \sum_{n \in \mathbb{N}} \mathbb{P}[\, e_n = 1 \,] \mathbb{P}_0[\, \tau_1 \wedge \tau_{-1} \in J_n \,] = \alpha.$$

Thus, $\gamma = \alpha$, and $X$ is the Skew Brownian motion of parameter $\alpha$. □

**Corollary 3.** *For all $t > 0$, $\mathbb{P}_0[\, X_t \geq 0 \,] = \alpha$.*

Let $\tau$ be the right continuous inverse of the local time $L_t^0(X)$ of a continuous diffusion process $X$. Each jump of $\tau$ corresponds to an excursion of this process, that is for each $t$ such that $\tau(t) = \tau(t-)$, there exists an interval $J_n = (\ell_n, r_n)$ such that $\ell_n = \tau(t-)$ and $r_n = \tau(t)$. Let $\mathcal{U}$ be the set of excursions, that is the set of continuous functions $f$ from $\mathbb{R}_+$ to $\mathbb{R}$ such that $f(0) = 0$, $f(t) = 0$ for all $t \geq \zeta$ for some $\zeta > 0$ and either $f(t) > 0$ for $t \in (0, \zeta)$ or $f(t) < 0$ for $t \in (0, \zeta)$.

For any $n \in \mathbb{N}$, $f(t) = X((t - r_n) \wedge (\ell_n - r_n))$ is some element of $\mathcal{U}$, with the life-time $\zeta = \ell_n - r_n$. Accordingly, for each $t \in \mathcal{J}$ where $\mathcal{J} = \{ t \geq 0 \mid \tau(t) \neq \tau(t-) \}$, one may associate a point $f(t)$ in $\mathcal{U}$, corresponding to some excursion of $X$. A striking result due to K. Itô is that the process $(t, f(t))_{t \geq 0}$ is a homogeneous Poisson point process with values on $\mathbb{R}_+ \times \mathcal{U}$ and intensity measure $\mathrm{d}t \times \widehat{\mathbb{P}}$. The measure $\widehat{\mathbb{P}}$ is a $\sigma$-finite measure, but which is infinite. This measure $\widehat{\mathbb{P}}$, called the excursions' measures, fully characterizes the diffusion process.

For a measurable subset $\Gamma$ of $\mathcal{U}$ such that $\widehat{\mathbb{P}}[\Gamma]$ is finite, $\widehat{\mathbb{P}}[\Gamma]$ denotes the average number of excursions in $\Gamma$ per unit of local time.

We denote by $\mathcal{U}^+$ (resp. $\mathcal{U}^-$) the set of positive (resp. negative) excursions, *i.e.* the excursions $f$ for which $f(t) > 0$ (resp. $f(t) < 0$) on $(0, \zeta)$.

On the excursions theory, see [11, 74]...

The following proposition is a direct result from the construction of the Skew Brownian motion given by the theory of excursions (see also [11]).

**Proposition 11.** *For the Skew Brownian motion of parameter $\alpha$, the measure $\widehat{\mathbb{P}}$ may be written*

$$\widehat{\mathbb{P}} = \alpha \widehat{\mathbb{P}}^+ + (1 - \alpha)\widehat{\mathbb{P}}^-, \tag{41}$$

*where $\widehat{\mathbb{P}}^+$ (resp. $\widehat{\mathbb{P}}^-$) is the excursions' measure of the Reflected Brownian motion on $\mathbb{R}_+$ (resp. on $\mathbb{R}_-$), that is $(|W_t|)_{t \geq 0}$, (resp. $(-|W_t|)_{t \geq 0}$), where $W$ is a Brownian motion.*

*Remark 8.* In [89] (see also [43]), S. Watanabe shows how to construct a diffusion given its excursions' measures, which provides another way to construct the SBM from Proposition 11.

*Remark 9.* This proposition has to be connected with the relations between the left, right and symmetric local times in (39).

*Remark 10.* For a diffusion $Y$ with a scale function $S^Y$, The function $1/S^Y(x)$ is also the average number of excursions reaching the point $x$ per unit of local time. Thus, this proposition is also coherent with the fact that if $S^B(x) = x$ is the scale function of the Brownian motion, the average number of excursions of the Skew Brownian motion reaching $x > 0$ (resp. $x < 0$) per unit of local time is equal to $\alpha/x$ (resp. $(1 - \alpha)/x$), and thus the scale function of the Skew Brownian motion is $S(x)$ given by (26).

*Remark 11.* The article [55] studies the link between the coefficients in a decomposition of type (41) for a general diffusion with coefficients discontinuous at 0, that is related to the derivatives of some Green functions at 0. We recover the decomposition (41) with $L = \frac{1}{2a(x)} \frac{\mathrm{d}}{\mathrm{d}x} \left( a(x) \frac{\mathrm{d}}{\mathrm{d}x} \right)$ with $a(x) = \alpha$ if $x \geq 0$ and $a(x) = 1 - \alpha$ if $x < 0$.

**Another construction of the density transition function.** The reflection principle for the Brownian motion together with the decomposition of excursions allows us to give another construction of density transition function of the Skew Brownian motion, following a result of J. Walsh [87]. Let $\tau$ be the first time the

Skew Brownian motion hit 0, and $B$ be a Brownian motion. Then

$$\mathbb{P}_x[\, X_t \in \mathrm{d}y\,] = \mathbb{P}_x[\, X_t \in \mathrm{d}y; \tau \le t\,] + \mathbb{P}_x[\, X_t \in \mathrm{d}y; \tau > t\,].$$

If $(x, y)$ belongs to $\mathbb{R}_+^* \times \mathbb{R}_-$ or to $\mathbb{R}_-^* \times \mathbb{R}_+$, by the continuity of the path, $\mathbb{P}_x[\, X_t \in \mathrm{d}y; \tau > t\,] = 0$. If $(x, y)$ belongs to $\mathbb{R}_+^* \times \mathbb{R}_+$ or to $\mathbb{R}_-^* \times \mathbb{R}_-$, then

$$\mathbb{P}_x[\, X_t \in \mathrm{d}y; \tau > t\,] = \mathbb{P}_x[\, B_t \in \mathrm{d}y; \tau > t\,] = (p(t, x, y) - p(t, x, -y))\,\mathrm{d}y,$$

where $p(t, x, y)$ is the transition density function of the Brownian motion. On the other hand, by the reflection principle and the construction of the Skew Brownian motion,

$$\mathbb{P}_x[\, X_t \in \mathrm{d}y; \tau \le t\,] = \begin{cases} \alpha \mathbb{P}_x[\, B_t \in -\mathrm{d}y\,] & \text{if } x > 0, \ y \ge 0, \\ (1 - \alpha)\mathbb{P}_x[\, B_t \in \mathrm{d}y\,] & \text{if } x \ge 0, \ y < 0, \\ \alpha \mathbb{P}_x[\, B_t \in \mathrm{d}y\,] & \text{if } x \le 0, \ y \ge 0, \\ (1 - \alpha)\mathbb{P}_x[\, B_t \in -\mathrm{d}y\,] & \text{if } x < 0, \ y < 0. \end{cases} \tag{42}$$

Hence, after a short computation, this gives us Formula (17).

## 8. Approximation by random walks (I)

As one can expect, a Skew Brownian motion may be approximated by the following random walk: Let $(S_k)_{k \ge 0}$ be the random walk starting from $0$ with probability transition

$$\begin{cases} \mathbb{P}[\, S_{k+1} = i + 1 \mid S_k = i\,] = \mathbb{P}[\, S_{k+1} = i - 1 \mid S_k = i\,] = \frac{1}{2} \text{ if } i \ne 0, \\ \mathbb{P}[\, S_{k+1} = 1 \mid S_k = 0\,] = \alpha, \\ \mathbb{P}[\, S_{k-1} = -1 \mid S_k = 0\,] = 1 - \alpha. \end{cases}$$

We set for all $t \ge 0$ and any integer $n$,

$$X_t^n = \frac{1}{n} S_{\lfloor n^2 t \rfloor} + \frac{n^2 t - \lfloor n^2 t \rfloor}{n} (S_{1 + \lfloor n^2 t \rfloor} - S_{\lfloor n^2 t \rfloor}).$$

The following theorem may be found in [41], and then in [52] and in [17] in a more general setting. P. Étoré used this result in [27] to construct a numerical scheme and compute its speed of convergence (see Section 11.9.3). We give another proof of this theorem in Section 9 by constructing the random walk $(S_n^k)_{k \in \mathbb{N}}$ from the trajectories of the SBM($\alpha$).

**Theorem 7** ([41]). *The sequence* $(X^n)_{n \in \mathbb{N}}$ *converges in distribution in the space of continuous functions to the SBM($\alpha$)* $X$.

*Proof.* Let us prove first the convergence of the marginals of $X^n$. For this, we remark that

$$\mathbb{P}[\, S_k = \ell\,] = \begin{cases} \alpha \mathbb{P}[\, |S_k| = \ell\,] & \text{if } \ell > 0 \\ (1 - \alpha)\mathbb{P}[\, |S_k| = -\ell\,] & \text{if } \ell < 0. \end{cases}$$

Hence,

$$\mathbb{E}[\exp(i\lambda X_t^n)] = \alpha\psi^n(\lambda) + (1-\alpha)\psi^n(-\lambda) \text{ with } \psi(\lambda) = \mathbb{E}[\exp(i\lambda|X_t^n|)].$$

Using the results of Section 7, one has also that

$$\mathbb{E}[\exp(i\lambda X_t)] = \alpha\psi(\lambda) + (1-\alpha)\psi(-\lambda) \text{ with } \psi(\lambda) = \mathbb{E}[\exp(i\lambda|X_t|)].$$

The results follows from the convergence of the normalized reflected random walk $|X^n|$ to the reflected Brownian motion $|X|$ given by the Donsker theorem.

The proof of the convergence finite-dimensional distributions of the $X^n$'s is similar and uses the same kind of computations as in (42). Yet the computations become heavy so that we skip it: in the next Section, we will see how to construct our random walk $(S_k^n)_{k\in\mathbb{N}}$ from a trajectory of $X$, and then how to prove the convergence of $X_t^n$ (which is then defined on the probability space of $X$) in probability to $X_t$. The convergence of the finite-dimensional distributions is then immediate.

The tightness of $(X^n)_{n\in\mathbb{N}}$ follows from the Kolmogorov criteria [46, Theorem 2.8, p. 53]. For this, let us note that

$$|X_t^n - X_s^n|^2 \leq \frac{3}{n^2}(n^2t - \lfloor n^2t\rfloor)^2(S_{1+\lfloor n^2t\rfloor} - S_{\lfloor n^2t\rfloor})^2$$
$$+ \frac{3}{n^2}(n^2s - \lfloor n^2s\rfloor)^2(S_{1+\lfloor n^2s\rfloor} - S_{\lfloor n^2s\rfloor})^2 + \frac{3}{n^2}(S_{\lfloor n^2t\rfloor} - S_{\lfloor n^2s\rfloor})^2 \quad (43)$$

Let us set $\xi_k = S_{k+1} - S_k$. Then all the $\xi_k$ are independent with variance 1, and $\mathbb{E}[\xi_k]$ is equal to $2\alpha - 1$ if $S_k = 0$ and to 0 if $S_k \neq 0$. Thus for any integer $p > q > 0$,

$$\mathbb{E}[(S_p - S_q)^2] = \sum_{k=q}^{p-1} \text{Var}(\xi_k) + \left(\sum_{k=q}^{p-1}\mathbb{E}[\xi_k]\right)^2$$
$$\leq p - q + (2\alpha - 1)(p-q)^2 \leq 2\alpha(p-q)^2.$$

Hence, it follows from (43) that for any $s, t \geq 0$,

$$\mathbb{E}[|X_t^n - X_s^n|^2] \leq 3\max\{2\alpha, 1\}(t-s)^2.$$

One deduces that for any $\gamma < 1/2$, there exists a random variable $K_\gamma^n$ such that

$$\sup_{0 \leq s < t \leq T}|X_t^n - X_s^n| \leq K_\gamma^n(t-s)^\gamma \text{ and } \sup_{n\in\mathbb{N}}\mathbb{E}[K_\gamma^n] < +\infty.$$

This implies the tightness of $(X^n)_{n\in\mathbb{N}}$. $\qquad\square$

*Remark* 12. It follows from this proof that the trajectories of the SBM($\alpha$) are $\gamma$-Hölder continuous for any $\gamma < 1/2$. Proposition 2.1 in [14] proves the following fact: almost surely, the modulus of continuity

$$w(X, \delta) = \sup_{s,t\in[0,T], \ |t-s|\leq\delta}|X_t - X_s|$$

of the strong solution $X$ on $[0, T]$ to (37) is smaller than $2w(B, \delta)$ for any $\delta > 0$.

In [17], it is proved that $(X_t^n, Y_t^n)_{t \geq 0}$ converges in distribution to $(B_t^\alpha,$ $\int_0^t f(B_s^\alpha) \, \mathrm{d}B_s^\alpha)_{t \geq 0}$, where $(X_t^n)_{t \geq 0}$ and $(Y_t^n)_{t \geq 0}$ are the linear interpolations of $X_{k/n}^n = (\sqrt{n})^{-1} S_{k/n}$ and $Y_{k/n}^n = \sum_{i=1}^k f(X_{(i-1)/n}^n)(X_{i/n}^n - X_{(i-1)/n}^n)$, and $f$ is a measurable, bounded function. This allows us to prove some results related to the "horizontal-vertical" random walk.

*Remark 13.* In [41], it is also noted that if given $S_i = 0$, the distribution of $S_{i+1}$ has the distribution of an integrable random variable $Z$ with values in $\mathbb{Z}$, then $n^{-1} S_{\lfloor n^2 t \rfloor}$ converges to a SBM($\alpha$) with $\alpha = \mathbb{E}[\, Z^+\,]/\mathbb{E}[\,|Z|\,]$.

## 9. Approximation by random walks (II)

We give another way to prove Theorem 7 that starts from a trajectory $X(\omega)$ of the SBM($\alpha$).

For the sake of simplicity, let us assume that $X_0 = 0$. Fix some integer $n$, and let us construct recursively from $X$ the following sequence of stopping times: $\tau_0^n = 0$ and

$$\tau_{k+1}^n = \inf \left\{ t > \tau_k^n \, \middle| \, X_t = \frac{\ell}{n} \text{ for some integer } \ell \neq \ell' \text{ will } X_{\tau_k^n} = \frac{\ell'}{n} \right\}.$$

Of course, as the trajectories of the SBM($\alpha$) are continuous, $|X_{\tau_{k+1}^n} - X_{\tau_k^n}| = n^{-1}$. In other words, the sequence $(\tau_k^n)_{k=0,1,\dots}$ records the successive passage times of $X(\omega)$ on the grid $\{\, k/n \,|\, k \in \mathbb{N} \,\}$. Finally, set $S_k^n = X_{\tau_k^n}$ and define $X^n$ by

$$X_t^n = \frac{1}{n} S_k^n + \frac{t - k\delta t}{n} (S_{k+1}^n - S_k^n) \text{ when } t \in [k\delta t, (k+1)\delta t]$$

where $\delta t = 1/n^2$.

With the scale function given in (26), $\mathbb{P}[\, X_{\tau_{k+1}^n} = (\ell+1)/n \,|\, X_{\tau_k^n}\,]$ is equal to $\alpha$ if $\ell = 0$ and to $1/2$ otherwise. With the strong Markov property of $X$, this means that the random walk $(S_k^n)_{k \in \mathbb{N}}$ is equal in distribution to the random walk previously constructed in Section 8.

This theorem is a "specialization" to the SBM($\alpha$) of Theorem 4.1 in [52].

**Theorem 8.** *The process $X^n$ constructed from $X$ as above converges uniformly on $[0, T]$ in probability to $X$ with respect to $\mathbb{P}$.*

*Proof.* We prove first the convergence in probability of $X_t^n$ to $X_t$ for any $t \in [0, T]$. Let us note first that since the reflected excursion of the SBM($\alpha$) are the same as the one of the Brownian motion, then for $k = 0, 1, \dots$, the increments $\tau_{k+1}^n - \tau_k^n$ are independent and have the same distribution (whatever the value of $\alpha$). In particular, $\mathbb{E}[\, \tau_{k+1}^n - \tau_k^n\,] = 1/n^2 = \delta t$ (this explain our choice of $\delta t$). The key relation is

$$X_{\lfloor n^2 t \rfloor/n^2}^n = \frac{1}{n} S_{\lfloor n^2 t \rfloor}^n = X_{\tau_{\lfloor n^2 t \rfloor}^n}.$$

Hence,

$$\mathbb{P}[\,|X_t^n - X_t| > C\,] \leq \mathbb{P}[\,|X_t^n - X_{\tau_{\lfloor n^2 t \rfloor}^n}^n| > C/2\,] + \mathbb{P}[\,|X_t - X_{\tau_{\lfloor n^2 t \rfloor}^n}| > C/2\,].$$

In the proof of Theorem 7, the tightness of $(X^n)_{n\in\mathbb{N}}$ has been established by showing it satisfies the Kolmogorov criteria. Hence, if $\mathrm{osc}(X^n,\delta)$ is the modulus of continuity of $X^n$, then for any $\delta > 0$ and any $\gamma < 1/2$,

$$
\mathbb{P}\big[\,|X_t^n - X_{\tau^n_{\lfloor n^2t\rfloor}}^n|\,>C\,\big] \leq \mathbb{P}\big[\,\mathrm{osc}(X^n,\delta) > C; |\tau^n_{\lfloor n^2t\rfloor} - t| < \delta\,\big]
$$

$$
+\,\mathbb{P}\big[\,\mathrm{osc}(X^n,\delta) > C; |\tau^n_{\lfloor n^2t\rfloor} - t| > \delta\,\big]
$$

$$
\leq \frac{\delta^\gamma}{C}\mathbb{E}\big[\,K^n_\gamma\,\big] + \mathbb{P}\big[\,|\tau^n_{\lfloor n^2t\rfloor} - t| > \delta\,\big], \quad (44)
$$

where $K^n_\gamma$ is the random $\gamma$-Hölder constant of $X^n$ which is known to satisfy $\sup_{n\in\mathbb{N}}\mathbb{E}\big[\,K^n_\gamma\,\big] < +\infty$. One may replace $X^n$ by $X$ in (44), and the result follows from the convergence in probability of $\tau^n_{\lfloor n^2t\rfloor}$ to $t$. By a scaling argument, $(\tau^n_{k+1} - \tau^n_k)_{k=0,\ldots,\lfloor n^2t\rfloor}$ is equal in distribution to $(n^{-2}\sigma_k)_{k=0,\ldots,\lfloor n^2t\rfloor}$, where the $\sigma^n_k$'s are independent, identically distributed with the distribution of the first exit time from $[-1,1]$ of the Brownian motion starting from 0. The $\sigma_k$ have a mean equal to 1 and a finite variance, so that the desired result is a consequence of the weak law of large numbers.

Now, to prove the uniform convergence in probability of $X^n$ to $X$, it is sufficient to note that for $0 = t_1 < t_2 < \ldots < t_\ell = T$ with $t_{i+1} - t_i < \delta$ for $i = 1,\ldots,\ell-1$,

$$
\sup_{t\in[0,T]} |X_t^n - X_t| = \sup_{i=1,\ldots,\ell-1}\ \sup_{t\in[t_i,t_{i+1}]} |X_t^n - X_t|
$$

$$
\leq \sup_{i=1,\ldots,\ell-1}\left(\sup_{t\in[t_i,t_{i+1}]}|X_t^n - X_{t_i}^n| + |X_{t_i}^n - X_{t_i}| + \sup_{t\in[t_i,t_{i+1}]}|X_{t_i} - X_t^n|\right)
$$

$$
\leq \mathrm{osc}(X^n,\delta) + \mathrm{osc}(X,\delta) + \sup_{i=1,\ldots,\ell-1}|X_{t_i}^n - X_{t_i}|.
$$

Hence, one deduces that $\lim_{n\to\infty}\mathbb{P}\big[\sup_{t\in[0,T]}|X_t^n - X_t| \geq C\,\big] = 0$ from the stochastic equicontinuity of $(X^n)_{n\in\mathbb{N}}$ and of $X$ by choosing $\delta$ (and then the points $\{t_i\}_{i=0}^\ell$ small enough, and afterwards using the convergence in probability of the $X_{t_i}^n$ to $X_{t_i}$ for $i = 1,\ldots,\ell$. $\qquad\square$

In order to set up a simulation scheme, the rate of convergence is computed in [27, 28], in the more general context of diffusion processes with discontinuous coefficients (besides, using a proper time increment, we are not bound to use a random walk on a regular grid as shown in [28]. See also points **D.** and **E.** in Section 11.9.3).

**Proposition 12** ([27, 28]). *For any $\gamma < 1/2$ and any $T \geq 0$, there exists a constant $C_{\gamma,T}$ such that $\sup_{t\in[0,T]}\mathbb{E}\big[\,|X_t^n - X_t|\,\big] \leq C_{\gamma,T}n^\gamma$.*

## 10. A "follow the leader" construction

We present here another construction based on techniques coming from the theory of continuous multi-armed bandits (See [60] for example). Indeed, we present here in a simpler case the proof given in [2] that deals to a variably skewed Brownian motion (See Section 11.6). Some of the results given in [17] are also related to this construction.

Let us set $\alpha \in (0,1)$, $\beta = 2\alpha - 1 \in (-1,1)$ and $\gamma = (1+\beta)/(1-\beta) > 0$.

The idea is to construct a process $(Z_t^1, Z_t^2, U_t^1)_{t \geq 0}$ with values in $(\mathbb{R}_+)^2 \times (\mathbb{R}_+)^2 \times \mathbb{R}_+$ from which one can get the positive excursions of the SBM($\alpha$) from $Z^1$, its negative excursions from $Z^2$ and its local time from $U$.

Indeed, $Z = (Z^1, Z^2)$ moves either horizontally or vertically when away from the line $\Upsilon = \{(x,y) \mid y = \gamma x\}$. The horizontal (resp. vertical) displacements of $Z$ away from $(U_t^1, \gamma U_t^1)$ correspond to a positive (resp. negative) excursions of the SBM($\alpha$) arising when its local times takes the value $L_t = (1+\gamma)U_t^1$. So, we can read "graphically" the signs and the heights of the excursions on the local time scale: See Figure 1. It is then possible to reconstruct a SBM($\alpha$) from $(Z^1, Z^2)$.

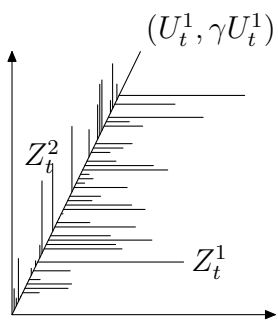

Fig 1. *A representation of the excursions of the SBM($\alpha$)*

The heuristic idea is to construct $Z^1$ and $Z^2$ as two time-changed independent Brownian motion $(B_{T^1(t)}^1)_{t \geq 0}$ and $(B_{T^2(t)}^2)_{t \geq 0}$ where either $T^1(t)$ or $T^2(t)$ increases at rate 1 while the other remains constant. The times at which the switches between $T^1$ and $T^2$ occur are obtained by comparing the supremum, as a function of time, of $\gamma B^1$ with the ones of $B^2$

We now turn to the rigorous construction.

Let us define a subset $D$ of $(\mathbb{R}_+)^2$ by

$$D = \left\{ (s_1, s_2) \in \mathbb{R}_+^2 \,\middle|\, \gamma \sup_{u_1 \leq s_1} B^1(u_1) > \sup_{u_2 \leq s_2} B^2(u_2) \right\},$$

whose closure $\overline{D}$ has the following properties: (a) $\mathbb{R}_+ \times \{0\} \subset D$; (b) $\{s \in \overline{D}\} \in \mathcal{F}_s$, where $(\mathcal{F}_{(s,t)})_{(s,t) \in \mathbb{R}_+^2}$ is the multi-parameter filtration defined by $\mathcal{F}_{(s,t)} = \mathcal{F}_s^1 \vee \mathcal{F}_t^2$; (c) If $(s_1, s_2) \in D$, then

$$\left\{ (u_1, u_2) \in \mathbb{R}_+^2 \,\middle|\, u_1 \geq s_1, \ 0 \leq u_2 \leq s_2 \right\} \subset \overline{D}.$$

It follows from these properties that $\overline{D}$ is described by the part of its boundary that intersects $(\mathbb{R}_+^*)^2$: See Figure 2. In addition, there exists a unique process $T(t) = (T^1(t), T^2(t))$ — called a *strategy* — such that [88] (a) $T(t) = (0, 0)$ and $T^1, T^2$ are non-decreasing; (b) $T^1(t) + T^2(t) = t$ for $t \geq 0$; (c) for all $s_1, s_2 \geq 0$, the set $\left\{ (s_1, s_2) \,\middle|\, T^1(t) \leq s_1, \ T^2(t) \leq s_2 \right\}$ is $\mathcal{F}_{(s_1, s_2)}$-measurable; (d) the graph of $(T(t))_{t \geq 0}$ parametrizes the boundary of $\overline{D}$.

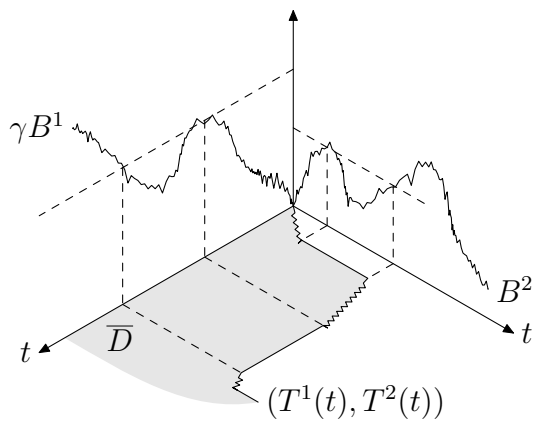

Fig 2. *Construction of the time changes*

Indeed, $T^1(t)$ increases at rate 1 when $B^2_{T^2(t)} > \gamma B^1_{T^1(t)}$ and $T^2(t)$ increases at rate 1 when $B^2_{T^2(t)} < \gamma B^1_{T^1(t)}$. In other words, the strategy "follows the leader" between $\gamma B^1$ and $B^2$.

Let us now consider a maximal time interval $(t, t')$ with $t' > t$ such that $T^1$ is strictly increasing on it, which means that $T^2$ is constant on $(t, t')$. We first note that as $(T^1(t), T^2(t))$ belongs to the boundary of $\overline{D}$,

$$\text{for any } u \in [t, t'], \ \gamma \sup_{s_1 \in [0, T^1(u)]} B^1_{s_1} = \sup_{s_2 \in [0, T^2(t)]} B^2_{s_2}.$$

As $(t, t')$ is chosen to be maximal, for all $u \in [t, t']$, $\sup_{s_1 \in [0, T^1(u)]} B^1_{s_1} = B^1_{T^1(t')} = B^1_{T^1(t)}$, which means that $(B^1_{T^1(u)})_{u \in [t, t']}$ performs an excursion below the level $B^1_{T^1(t)}$ on the $[t, t']$. As $(T^1(t), T^2(t))$ belongs to the boundary of $\overline{D}$, there exists no time $\widehat{t} < T^2(t)$ such that $\sup_{s \in [0, \widehat{t}]} B^2_s = B^1_{T^1(t)}$. This means that $B^2_{T^2(t)} = \sup_{t \in [0, T^2(t)]} B^2_s$. Hence, at time time $t$, $B^1_{T^1(t)}$ has a local maximum and $B^2_{T^2(t)}$ reaches its maximum.

We now set for $r \geq 0$, $Z^i_r = B^i_{T^i(r)}$ for $i = 1, 2$ and $U^1_r = \sup_{s \in [0, s]} Z^1_s$, $U^2_r = \gamma U^1_r$. From the previous considerations,

$$\gamma U^2_t = \gamma U^1_{t'} = \gamma Z^1_t = \gamma Z^1_{t'} = Z^2_t = Z^2_{t'} = U^2_t = U^2_{t'}.$$

Thus, on the time interval $[t, t']$ above, $Z^1$ performs an excursion below the level $U^1_t$ and then $Z$ moves horizontally away from $\Upsilon$, with $Z_t = Z_{t'} = U_t \in \Upsilon$.

A similar study could be done for $T^2$, and we see now how to define our candidate to be a SBM:

$$X_t = (\gamma U_t^1 - Z_t^2) - (U_t^1 - Z_t^1) \text{ and } L_t = (1+\gamma)U_t^1. \tag{45}$$

**Proposition 13.** *The process $X$ is a SBM($\alpha$) and $L$ is its local time.*

*Proof.* Let $(\mathcal{G}_t)_{t\geq 0}$ be the filtration defined by $\mathcal{G}_t = \mathcal{F}_{T(t)}$. Let us remark first that $B_t = -Z_t^2 + Z_t^1$ has for quadratic variation $\langle B \rangle_t = T^1(t) + T^2(t) = t$ from the very definition of the strategy and is $(\mathcal{G}_t)_{t\geq 0}$-adapted. By the Paul Lévy's theorem, $B$ is a Brownian motion and then $X_t = B_t + \beta L_t$ with $\beta = 2\alpha - 1$. The process $|X|$ is equal to the distance of $Z$ from $\Upsilon$ and then, by construction,

$$|X_t| = \widehat{B}_t + L_t$$

with $\widehat{B}_t = -B_{T^1(t)}^1 - B_{T^2(t)}^2$. Again with computing the quadratic variation of this latter process, $\widehat{B}$ is a $(\mathcal{G}_t)_{t\geq 0}$-Brownian motion. Besides,

$$-\widehat{B}_t = Z_t^1 + Z_t^2 = \int_0^t \operatorname{sgn}(X_s)\,\mathrm{d}B_s = \int_0^t \operatorname{sgn}(X_s)\,\mathrm{d}X_s$$

where $\operatorname{sgn}(x) = 1$ if $x > 0$, $\operatorname{sgn}(x) = -1$ if $x < 0$ and $\operatorname{sgn}(x) = 0$ if $x = 0$, as $\int_0^t \mathbf{1}_{\{Z_s \notin \Upsilon\}}\mathrm{d}L_s = 0$ if $X_t \neq 0$. Thus, with the Itô-Tanaka formula, $L_t$ is the symmetric local time of $X$. Finally, from (45), $X_t = B_t + (\gamma-1)/(\gamma+1)L_t$, and the skewness parameter $\alpha$ is identified from value of $\gamma$. $\qquad\square$

## 11. Applications and extensions

### 11.1. General constructions of "Skew" diffusion processes

The general theory of one-dimensional diffusion processes gives all the tools to construct rather easily a "skew" diffusion process. There are basically two ways to do so.

First, consider diffusion processes living on $[0, +\infty)$, or symmetric processes, and the sign of each excursion is changed with independent Bernoulli random variables. For this, one may use the description of a strong Markov process in term of a *minimal process* (*i.e.*, a process that is killed when it reaches 0) and an *entrance law* (See [11] for a whole account on this theory). See [86] for related results and constructions.

In [4, 5] a Skew Bessel process is constructed this way.

Indeed, the symmetry of the process plays an important role in the previous construction, so here is another one. Let $X$ be a process with a scale function $S$ and a speed measure $m$, $x_0$ in the state space of $X$ and $\gamma > 0$ let us define a new process $X^\gamma$ by the scale measure $S^\gamma$ and speed measure $m^\gamma$ defined by

$$S^\gamma(x) = \begin{cases} \gamma S(x) & \text{if } x \geq x_0, \\ S(x) & \text{if } x < x_0, \end{cases} \text{ and } m^\gamma(\mathrm{d}x) = \begin{cases} \gamma^{-1}m(\mathrm{d}x) & \text{if } x \geq x_0, \\ m(\mathrm{d}x) & \text{if } x < x_0 \end{cases}$$

if $S$ is such that $S(x_0) = 0$. As the behavior of a process around a point $x$ is specified locally by its scale function and speed measure and the process generated by $(S, m)$ is the same as the process generated by $(\lambda S, \lambda^{-1} m)$ for all $\lambda > 0$, the process $X^\gamma$ behaves like $X$ except at $x_0$. The "skewness" of this process may be appreciated by the result stated in Remark 10.

The article [21] shows how to construct an *asymmetric Bessel process* by this method and provides applications to the computation of Asian options with discontinuous volatility. In [22], various processes are constructed by this way, in order to model some prices whose volatility suddenly increases or decreases when the price goes above or below a threshold.

### 11.2. Walsh's Brownian motion

Initially introduced by J. Walsh in [87], the *Walsh's Brownian motion* is a diffusion process on a set of $n$ rays in $\mathbb{R}^2$ emanating from 0. To each ray $I_i$ is associated a weight $\alpha_i$ corresponding heuristically to the probability for the process to go in this ray. On each ray, the process behaves like a Brownian motion.

Of course, due to the irregularities of the trajectories of the Brownian motion, this description is a non-sense, but this process may be described by its excursions' measure

$$\widehat{\mathbb{P}} = \sum_{i=1}^{n} \alpha_i \widehat{\mathbb{P}}^i, \text{ with } \sum_{i=1}^{n} \alpha_i = 1 \text{ and } \alpha_i \in (0, 1),$$

where $\widehat{\mathbb{P}}^i$ is the excursion measure of the Reflected Brownian motion on the ray $I_i$. On this process, see also [4].

The Skew Brownian motion corresponds to the case $n = 2$. The *spider martingale* generalizes the notion of Walsh's Brownian motion. This object has given rise to an abundant literature on Brownian filtrations, especially by giving a negative answer to the question "if a Brownian motion is adapted to some filtration, is this filtration generated by a Brownian motion?" (See [7, 85], [94, Sect. 17, p. 103] and the cited works within), but this goes far beyond the scope of this article.

### 11.3. A generalized arcsine law

For the Walsh's Brownian motion of Section 11.2, let $A^i(t)$ be the occupation time up to time $t$ of the ray $I_i$, that is $A^i(t) = \int_0^t \mathbf{1}_{\{X(s) \in I_i\}} \, \mathrm{d}s$. As for the Brownian motion, $A^i(t) \overset{\text{dist.}}{=} t A^i(1)$ for any $t \geq 0$. The arcsine law for the Brownian motion may be extended to this case.

**Theorem 9** ([5]). *Let $T_1, \ldots, T_n$ be independent, stable random variables of index $1/2$. Then*

$$(A_1(1), \ldots, A_n(1)) \stackrel{\text{dist.}}{=} \left( \frac{\alpha_i T_i}{\sum_{j=1}^n \alpha_j T_j} \right)_{i=1,\ldots,n}.$$

More is proved in [5], and Skew Bessel processes are also considered.

Assume that $n = 2$ and that $T_1$, $T_2$ are two independent copies of a stable random variable with exponent $\gamma \in (0, 1)$ (*i.e.*, $\mathbb{E}[\exp(-\zeta T_i)] = \exp(-\zeta^\gamma)$ for any $\zeta > 0$). Define by $\xi_{\gamma,\alpha}$ the random variables

$$\xi_{\gamma,\alpha} = \frac{\alpha^{1/\gamma} T_1}{\alpha^{1/\gamma} T_1 + (1-\alpha)^{1/\gamma} T_2}.$$

The distribution function of $\xi_{1/2,\alpha}$ is

$$\mathbb{P}[\, \xi_{1/2,\alpha} \leq x\,] = \frac{2}{\pi} \arcsin \sqrt{\frac{x}{x + (\alpha/(1-\alpha))^2 (1-x)}}.$$

For $\alpha = 1/2$, this is the arcsine distribution.

In Theorem 9, the occupation time of $\mathbb{R}_+$ or $\mathbb{R}_-$ at time 1 of the Skew Brownian motion is a random variable of type $\xi_{1/2,\alpha}$.

If the Skew Brownian motion is replaced by a Skew Bessel process of dimension $2 - 2\nu$, then the occupation time at time 1 is distributed as one of the $\xi_{\nu,\alpha}$ [5, 90]. In [90], it is also proved that for a general one-dimensional diffusion, the possible limit of $A_{\mathbb{R}^+}(t)/t$, where $A_{\mathbb{R}^+}(t)$ is the occupation time of $\mathbb{R}^+$, is a random variable of type $\xi_{\nu,\gamma}$.

In the case of a Brownian motion on a graph (See below Section 11.9.2), analytical considerations allows us also to compute the Laplace transforms of the occupation time on each of the edges of the graph [23].

### *11.4. The Skew Brownian motion as a flow*

As the SBM is a strong solution to the SDE,

$$X_t^z = z + B_t + \beta L_t^0(X^z) \tag{46}$$

it is natural to study the dependence of the SBM($\alpha$) with respect to its starting point.

**Proposition 14** ([3]). *For $x, y \in \mathbb{R}$, let $X^x$ and $X^y$ be the two strong solutions of (46) with the same parameter $\beta \in [-1, 1] \backslash \{0\}$ and the same Brownian motion $B$. Then almost surely, there exists a finite time $\tau$ such that $X_\tau^x = X_\tau^y$.*

From this coupling property, we deduce the following: for $\beta \in (-1, 1) \setminus \{0\}$ and $x \neq y \in \mathbb{R}$, there exists almost surely some $\tau > 0$ such that $x + \beta L_\tau(X^x) = y + \beta L_\tau(X^y)$.

If $K_t^x(X^x) = x + \beta L_t(X^x)$, the uniqueness of the strong solution to (47) implies that $K_s^x(X^x) = K_s^y(X^x)$ for $s \geq t$.

From this property, K. Burdzy and Z.-Q. Chen give in [13] several results related to the distributions and behavior of the processes $x \in \mathbb{R} \mapsto L_\tau^0(X^x)$ with $\tau = \inf \left\{ t \geq 0 \,\big|\, L_t^0(X^0) = 1 \right\}$ and $\beta \in \mathbb{Q} \cap (-1, 1) \mapsto \beta L_t^0(X^0)$.

It is also possible to extend (46) to $t \in \mathbb{R}$ by defining the Brownian motion $(B_t)_{t \in \mathbb{N}}$ as $B_t = \mathbf{1}_{\mathbb{R}_+}(t)W_t + \mathbf{1}_{\mathbb{R}_-^*}(t)W'_{-t}$ for two independent Brownian motions $W$ and $W'$. Let us note that $t \mapsto L_t^0(X^x)$ is a non-decreasing function on $\mathbb{R}$. If $X^x$ is the solution to (46), then $Y_t^x = X_{-t}^x$ is solution to (46) with $B_t$ replaced by $B_{-t}$ and $\beta$ by $-\beta$. Then $K_t^x(Y^x)$ also satisfies the coupling property on $\mathbb{R}_+$. This implies that almost surely, there exists some time $\tau < 0$ such that for $t \leq \tau$, $K_t^x(X^x) = K_t^y(X^y)$ and $K_t^x(X^x) \neq K_t^y(X^y)$ for $t > \tau$.

This phenomena is studied in [13, 14], where the authors consider the following generalization of (46), which is for $s, x \in \mathbb{R}$ the solution $(X^{s,x}, \ell^{s,x})$ to

$$
\begin{cases}
X_t^{s,x} = x + B_t - B_s + \beta \ell_t^{s,x}, \ t \geq s \\
\text{with } \ell_t^{s,x} = \lim_{\varepsilon \to 0} \frac{1}{2\epsilon} \int_0^t \mathbf{1}_{\{X_s \in [-\varepsilon, \epsilon]\}} \, \mathrm{d}s,
\end{cases}
\tag{47}
$$

Among many results, it is shown in [14] that (47) may have for some $(s, x) \in \mathbb{R}^2$ three distinct strong solutions. The notion of "lenses" allows us to characterize the points $(s, x)$ for which the decoupling/coupling of $K_t^x(X^x)$ takes place.

It is then conjectured that, although there exists almost surely a strong solution to $X_t^x = x + B_t + \beta L_t^0(X^x)$ simultaneously for all $x \in \mathbb{Q}$, there is no strong uniqueness for a solution defined simultaneously for all $x \in \mathbb{R}$. We refer the readers to the articles cited above for a complete account on these developments, results and proofs.

### 11.5. Comparison principle

A natural question is to know whether or not a comparison principle holds for the SBM. We give first two propositions. The first one concerns the case where the skewness parameters $(1 + \beta)/2$ are the same, but gives a strong statement on the solutions to (47).

**Proposition 15** ([13, Proposition 1.8]). *For all $(s, s', x, x') \in \mathbb{Q}^4$ simultaneously with a fixed $\beta \in [-1, 1] \setminus \{0\}$, either $X_t^{s,x} \leq X_t^{s',x'}$, $t \in \mathbb{R}$ or $X_t^{s,x} \geq X_t^{s',x'}$, $t \in \mathbb{R}$ almost surely, where $X^{s,x}$ is solution to (47).*

In particular, if $x \leq y$, then $X_t^x \leq X_t^y$ $t \geq 0$, where $X^x$ is the strong solution to (46). The second proposition concerns the case of different parameters.

**Proposition 16** ([52, Theorem 3.4],[96]). *Let $X^1$ and $X^2$ be two strong solutions to (46) respectively with starting $x^1$ and $x^2$ and parameters $\beta_1$ and $\beta_2$. If $x^1 \leq x^2$ and $\beta_1 \leq \beta_2$, then $X_t^1 \leq X_t^2$ for all $t > 0$.*

The distance between two solutions may then be computed.

**Theorem 10** ([96]). *Let $X^1$ and $X^2$ be two strong solutions to (46) with either:*
*(a) The same starting point $x$ and two different parameters $\beta_1$ and $\beta_2$ in $(-1, 1)$.*
*(b) Different starting points $x^1$ and $x^2$ and the same parameter $\beta \in (-1, 1)\backslash\{0\}$.*
*Then, the $\mathrm{L}^1$ distance between $X^1$ and $X^2$ is given by either*

$$(a) \ \mathbb{E}[\,|X_t^1 - X_t^2|\,] \leq |\beta_1 - \beta_2| \int_0^t \frac{1}{\sqrt{2\pi s}} I(t, x) \text{ with } I(t, x) = \exp\left(-\frac{x^2}{2s}\right) \mathrm{d}s,$$

$$(b) \ \begin{cases} \mathbb{E}[\,|X_t^1 - X_t^2|\,] \leq |x^1 - x^2| + |\beta| \cdot |I(t, x^1) - I(t, x^2)|, \\ \mathbb{E}[\,|L_t^0(X^1) - L_t^0(X^2)|\,] \leq \frac{1}{|\beta|}|x^1 - x^2| + |I(t, x^1) - I(t, x^2)|. \end{cases}$$

### 11.6. Generalization of the Skew Brownian motion (I)

A possible generalization of the Skew Brownian motion is to use a variable "skewness coefficients", which leads to an SDEs with local time of type

$$X_t = x + B_t + \beta(L_t^0(X)), \tag{48}$$

where $\beta$ is a differentiable function on $\mathbb{R}_+$ with $\beta'(x) \in (-1, 1)$ for all $x \in \mathbb{R}_+$. This SDE has been studied in [2].

**Theorem 11** ([2]). *There exists a unique strong solution to (48).*

The proof of the theorem relies on an extension of the construction given in Section 10. The variably skew Brownian motion has been used recently for the Skorokhod embedding problem [18].

### 11.7. Generalization of the Skew Brownian motion (II)

Another natural generalization is the solution of the SDE

$$X_t = x + \int_0^t \sigma(X_s)\,\mathrm{d}B_s + \int_0^t b(X_s)\,\mathrm{d}s + \int_0^t \alpha(s)\,\mathrm{d}L_s^0(X), \tag{49}$$

which is motivated by some homogenization — *i.e.*, change of scale — results [92].

**Theorem 12** ([91]). *If $\sigma$ is bounded below by a positive constant, bounded and Lipschitz continuous, $b$ is bounded and $\alpha$ is continuous with $|\alpha| \leq 1$, then the solution of (49) is pathwise unique.*

*Remark 14.* In [91], S. Weinryb uses a normalization of the local time different from ours.

### 11.8. Generalization of the Skew Brownian motion (III)

The SDE solved by the Skew Brownian motion may be generalized as a SDE of type

$$X_t = X_0 + \int_0^t \sigma(X_s)\,\mathrm{d}B_s + \int_{\mathbb{R}} \nu(\,\mathrm{d}x)L_t^x(X), \tag{50}$$

where $\nu$ is a signed measure. Of course, the case $\sigma = 1$ and $\nu = (2\alpha - 1)\delta_0$ corresponds to the SBM($\alpha$).

Let BV be the space of functions $f : \mathbb{R} \to \mathbb{R}$ of bounded variation such that i) $f$ is right continuous, ii) There exists $\varepsilon > 0$ such that $f(x) \geq \epsilon$.

Let M be the space of all signed measures $\nu$ on $\mathbb{R}$ such that $|\nu(\{x\})| < 1$ for all $x \in \mathbb{R}$.

**Theorem 13** ([9, 24, 51, 52]). *There exists a unique strong solution to* (50) *when* $\sigma \in$ BV *and* $\nu \in$ M.

Of course, the proof relies on an argument similar to the one given in Section 5.2 on the SDE that the SBM($\alpha$) solves. Here, the idea is to construct a harmonic function by the following way [52, Lemma 2.1]: for $\nu \in$ M, there exists a unique (up to a multiplicative constant) such that

$$f'(\,\mathrm{d}x) + (f(x) + f(x-))\nu(\,\mathrm{d}x) = 0$$

and, with the additional assumption that $f(x) \to 1$ when $x \to -\infty$,

$$f_\nu(x) \stackrel{\text{def.}}{=} f(x) = \exp(-2\nu^c((-\infty, x])) \prod_{y \leq x} \frac{1 - \nu(\{y\})}{1 + \nu(\{y\})}. \tag{51}$$

As $f$ is in BV, $f'(\,\mathrm{d}x)$ is the bounded measure associated to $f$ as in (28). Moreover, $\nu^c$ is the continuous part of $\nu$. In some sense, the measure $\nu$ corresponds to the logarithmic derivative of the function $f(x)$. Then, we set $F_\nu(x) = \int_0^x f_\nu(x)\,\mathrm{d}x$. It follows that $X$ is solution to (50) if and only if $Y = F_\nu(X)$ is solution to

$$Y_t = F_\nu(X_0) + \int_0^t (\sigma f_\nu)(F_\nu^{-1}(Y_s))\,\mathrm{d}B_s.$$

The existence of a weak solutions follows from the martingale problem [79, 46] while the pathwise uniqueness follows from a result due to S. Nakao [63]. The Yamada-Watanabe theorem allows us to conclude (See [46] for example).

For the SBM($\alpha$), $F_\nu(x)$ corresponds, up to a multiplicative function, to the scale function $S(x)$ given in (26).

In [52], J.-F. Le Gall also provides the following convergence theorem.

**Theorem 14** ([52, Theorem 3.1]). *Let* $(\nu_n)_{n\in\mathbb{N}}$ *be a sequence of measures in* M *and* $(\sigma_n)_{n\in\mathbb{N}}$ *be a sequence of functions in* BV *such that for some constants* $M$ *and* $\varepsilon$, i) $\sup_{n\in\mathbb{N}} |\nu_n(\mathbb{R})| \leq M$, ii) $\varepsilon \leq \sigma_n(x) \leq M$ *for all* $n \in \mathbb{N}$, $x \in \mathbb{R}$ *and* $\sup_{n\in\mathbb{N}} |\nu_n(\{x\})| \leq 1 - \varepsilon$ *for all* $x \in \mathbb{R}$. *Let us denote by* $X^n$ *the solution to* (37) *with* $\sigma$ *and* $\nu$ *replaced by* $\sigma_n$ *and* $\nu_n$, *with respect to the same Brownian motion (and then on the same probability space) and assume that* $X_0^n$ *converges in* $\mathrm{L}^1$ *to* $X_0$.

*Assume moreover, that there exist two functions* $\sigma$ *and* $f$ *in* BV *such that*

$$\sigma_n \xrightarrow[n\to\infty]{\mathrm{L}^1_{\mathrm{loc}}(\mathbb{R})} \sigma \text{ and } f_{\nu_n} \xrightarrow[n\to\infty]{\mathrm{L}^1_{\mathrm{loc}}(\mathbb{R})} f$$

*and set* $\nu(\mathrm{d}x) = -f'(\mathrm{d}x)/(f(x) + f(x-))$. *Then* $X^n$ *converges uniformly in* $\mathrm{L}^1(\mathbb{R})$ *to* $X$ *on* $[0, T]$ *for any* $T > 0$, *where* $X$ *is the strong solution to* $X_t = X_0 + \int_0^t \sigma(X_s)\,\mathrm{d}B_s + \int_{\mathbb{R}} \nu(\mathrm{d}x)L_t^x(X)$.

*Remark* 15. This theorem endows the importance of the function $f$ and allows us to recover the result of Proposition 8. It has to be noted that one may construct a sequence $(\nu_n)_{n \in \mathbb{N}}$ of measures in M satisfying the hypotheses of Theorem 14 converging weakly to a measure $\mu$ but for which the measure $\nu$ given by this theorem differs from $\mu$. Proposition 8 illustrates this points if we take for drift $b$ a symmetric function with compact support on $[-1, 1]$ with $\int_{-\infty}^{\infty} b(x) = \kappa/2$, so that $nb(x/n)\,\mathrm{d}x$ converges to the Dirac measure $(\kappa/2)\delta_0$. Yet $n^{-1}X_{tn^2}$ converges to the solution to $Y_t = B_t + \frac{e^\kappa - 1}{e^\kappa + 1}L_t^0(Y)$.

J.-F. Le Gall also proved some extension of the Donsker theorem when $\sigma = 1$ in (50). This Donsker theorem was used and improved by P. Étoré in [27, 28, 29] in order to approximate one-dimensional diffusion processes with discontinuous coefficients by random walks (See also Section 11.9.3).

Recently, R. Bass and Z. Q. Chen have shown in [9] strong existence and pathwise uniqueness of (50) under the assumptions that i) $\nu \in$ M, ii) $a$ is bounded below and above by a positive constant iii) there exists a strictly increasing function $f$ such that

$$|\sigma(x) - \sigma(y)|^2 \leq |f(x) - f(y)|, \ \forall x, y \in \mathbb{R}. \tag{52}$$

Let us note that if $\sigma \in$ BV, then $\sigma$ satisfies (52) with $f(x) = 2\|\sigma\|_\infty \int_{-\infty}^x |\,\mathrm{d}\sigma|$, so that the results of [9] are slightly more general than the one in [52].

The article [9] also provides us with a comparison principle and shows that there is no solution to (50) when $\nu(\{x\}) > 1$ for some $x \in \mathbb{R}$ (the non-existence of a solution was stated in [52] without proof).

As the local time plays a central role in the theory of one-dimensional diffusion processes, SDE with local times are also considered in the study of diffusions with *singular drift*, in order to get some diffusions under weak hypotheses (see for example [9, 10, 30], ...), Dirichlet processes [8, 25], ... In addition, they also appear as a tool to study the "critical cases" for results on strong/weak existence and uniqueness of solutions of SDE (see [24, 52] for example).

Finally, let us note that this kind of process appears also naturally when one studies diffusion processes generated by differential operators with discontinuous coefficients: see Section 11.9.

### 11.8.1. Spatial change of variable

The symmetric Itô-Tanaka formula allows us to state a stability result of solutions of (50) under a subclass of functions in BV.

**Proposition 17.** *Let* $F$ *be a function whose derivative* $f$ *belongs to* BV *and is bounded above and below by a positive constant constant. If* $X$ *is the solution to*

(50) *with $\sigma \in \mathrm{BV}$ and $\nu \in M$, then*

$$F(X_t) = F(X_0) + \int_0^t f(X_s)\sigma(X_s)\,\mathrm{d}B_s$$
$$+ \frac{1}{2}\int_{\mathbb{R}} ((f(x+) + f(x-))\nu(\,\mathrm{d}x) + f'(\,\mathrm{d}x))\,L_t^x(X). \quad (53)$$

*If $F$ is one-to-one, then $Y = F(X)$ is also solution to some SDE of type* (50) *with a coefficient $(f\sigma) \circ F^{-1} \in \mathrm{BV}$ and a measure $\mu \in \mathrm{M}$ with*

$$\mu(\,\mathrm{d}x) = \frac{g'(\,\mathrm{d}x)}{g(x+) + g(x-)} \text{ with } g = h \circ F^{-1} \text{ where } \nu(\,\mathrm{d}x) = \frac{h'(\,\mathrm{d}x)}{h(x) + h(x-)}$$

*for some function $h$.*

*Proof.* Formula (53) a direct consequence of the symmetric Itô-Tanaka formula (32). After a one-to-one change of variable, if $\varphi$ belongs to BV, then it is easily checked that $\int_{\mathbb{R}} \varphi'(\,\mathrm{d}x)k(x) = \int_{\mathbb{R}} \varphi'(\,\mathrm{d}x)k(F^{-1}(x))\,\mathrm{d}x$ for any bounded, measurable function $k$. The measure $\mu$ is then computed with Lemma 2.1 in [52] which asserts the existence of $h$, the previous change of variable and the relation $L_t^{F^{-1}(x)}(X) = L_t^x(F(X))$. $\qquad \square$

As we saw it, the existence and uniqueness results on SDEs with local time relies on some application of this formula, to get rid of the local time. In [6], M. Barlow used this kind of transform to construct counter-examples to the uniqueness of strong solutions to some SDEs.

From a numerical point of view, this change of variable is important in order to simulate the process $X$, as shown in Section 11.9.3.

In [67], Y. Ouknine made use of the Itô-Tanaka formula to study the distribution of solutions of SDEs of type (50) with $\sigma$ constant on $\mathbb{R}_+$ and on $\mathbb{R}_-^*$ and $\nu = \kappa\delta_0$. This class of SDEs remains stable when one uses for $F$ in Proposition 17 functions of type $F(x) = px^+ + qx^-$.

### 11.8.2. Random time change

Random time change is also a tool to study SDE with local time, as shown first in [24]. In [61, Chapter 5], M. Martinez used a random time change based on the SBM instead of a Brownian motion to provide us with existence and uniqueness results of one-dimensional SDEs in the spirit of the work of H. Engelbert and J. Schmidt (see [24] or [46, Section 5.5.A]).

## 11.9. Diffusions with discontinuous coefficients and diffusions on graphs

We present now in this section some results about PDEs with discontinuous coefficients, and diffusion equation on graphs. This gives the general framework to consider both a large subclass of SDEs of type (50) and extensions of the Walsh's Brownian motion presented in Section 11.2.

### 11.9.1. Infinitesimal generators and associated diffusion processes

Let $a$ and $\rho$ be two measurable functions on $\mathbb{R}$ such that $0 < \lambda \le a(x) \le \Lambda$ and $\lambda \le \rho(x) \le \Lambda$ for all $x \in \mathbb{R}$. Let also $b$ be a measurable bounded

Let $(L, \mathrm{Dom}(L))$ be the differential operator

$$L = \frac{\rho}{2}\frac{\mathrm{d}}{\mathrm{d}x}\left(a\frac{\mathrm{d}}{\mathrm{d}x}\right) + b\frac{\mathrm{d}}{\mathrm{d}x}, \ \ \mathrm{Dom}(L) = \left\{\, f \in \mathcal{C}(\mathbb{R};\mathbb{R}) \,\middle|\, Lf \in \mathcal{C}(\mathbb{R};\mathbb{R}) \,\right\}. \quad (54)$$

Using scale functions and speed measures, or analytical results on the fundamental solution of $\frac{\partial}{\partial t} - L$, or Dirichlet forms, it can be proved that $(L, \mathrm{Dom}(L))$ is the infinitesimal generator of a continuous, conservative, strong Markov process $(X_t, t \ge 0; \mathcal{F}_t, t \ge 0; \mathbb{P}_x, x \in \mathbb{R})$ (See [39, 54, 78] for example).

The next proposition follows from an application of Corollary 2 and Theorem 14, (See [27, 29, 54]). Note that a localization argument shall be used if $b$ does not belong to $\mathrm{L}^1(\mathbb{R};\mathbb{R})$ (see [78, Lemma II.1.12] for example).

**Proposition 18.** *We assume in addition to the previous hypotheses that $\rho$ and $\sigma$ belong to* BV. *Then for any starting point $x$, $X$ is the unique strong solution to the SDE*

$$X_t = x + \int_0^t \sqrt{a(X_s)\rho(X_s)}\,\mathrm{d}B_s + \int_0^t \frac{1}{2}a'(X_s)\rho(X_s)\,\mathrm{d}s$$
$$+ \int_0^t b(X_s)\,\mathrm{d}s + \int_{\mathbb{R}} \nu(\,\mathrm{d}y)L_t^y(X), \quad (55)$$

*where $B$ is a Brownian motion, $L_t^y(X)$ is the symmetric local time of $X$ at the point $y$ and*

$$\nu = \sum_{n \in J} \frac{a(x_n+) - a(x_n-)}{a(x_n+) - a(x_n-)}\delta_{x_n} \ \text{if} \ \left\{\, x_n \,\right\}_{n \in J \subset \mathbb{N}} = \left\{\, x \in \mathbb{R} \,\middle|\, a(x) \ne a(x-) \,\right\}.$$

The SDE (55) is a particular case of (50) encountered in Section 11.8.

The article [40] shows how to compute explicitly the density transition function and the Green function of this process using spectral analysis.

### 11.9.2. Transmission conditions and diffusions on a graph

The solution of the parabolic PDE $\frac{\partial u}{\partial t} = Lu$ with $u(0, x) = f(x)$, where $L$ is given by (54), has to be understood as a *weak solution*, that is as a function $u(t, x) \in \mathcal{C}(0, T; \mathrm{L}^2(\mathbb{R})) \cap \mathrm{L}^2(0, T; \mathrm{H}^1(\mathbb{R}))$ such that for all $\psi \in \mathcal{C}^\infty([0, T]; \mathbb{R})$,

$$\int_0^T \int_{\mathbb{R}} \frac{\partial \psi(t, x)}{\partial t} u(t, x)\rho(x)^{-1}\,\mathrm{d}x\,\mathrm{d}t - \int_0^T \int_{\mathbb{R}} b(x)\frac{\partial u(t, x)}{\partial x}\psi(t, x)\,\mathrm{d}x\,\mathrm{d}t$$
$$+ \int_0^T \int_{\mathbb{R}} a(x)\frac{\partial u(t, x)}{\partial x}\frac{\partial \psi(t, x)}{\partial x}\,\mathrm{d}x\,\mathrm{d}t = -\int_{\mathbb{R}} \varphi(x)\psi(0, x)\rho(x)^{-1}\,\mathrm{d}x.$$

Now, if the coefficients $a$ and $\rho$ belong to BV and the points $\{x_i\}_{i \in J}$ at which $a$ is discontinuous, then $u$ is also solution to the transmission problem

$$\begin{cases} \partial_t u(t,x) = Lu(t,x) \text{ on } \mathbb{R}_+^* \times (\mathbb{R} \setminus \{ x_i \}_{i \in J}), \\ u(0,x) = f(x) \text{ on } \mathbb{R}, \\ a(x_i-)\partial_x u(t,x_i-) = a(x_i+)\partial_x u(t,x_i+) \text{ for } i \in J, \\ u(t,x_i-) = u(t,x_i+) \text{ for } i \in J, \end{cases} \tag{56}$$

where $u \in \mathcal{C}^{1,2}(\mathbb{R}_+^* \times (\mathbb{R} \setminus \{ x_i \}_{i \in J})) \cap \mathcal{C}(\mathbb{R}_+^* \times \mathbb{R})$ (See [49]).

*Remark* 16. Note that we can get any arbitrary condition transmission at some point $x_i$ without changing the diffusion coefficient on the intervals $(x_{i-1}, x_i)$ and $(x_i, x_{i+1})$. For this, it is sufficient to multiply the coefficients $a$ by a $\mu_j$ and $\rho$ by $1/\mu_j$ on $(x_j, x_{j+1})$ for a well chosen sequence $\{ \mu_j \}_{j \in J}$.

Differential operators on graphs are natural expansions of one-dimensional differential operators with a transmission conditions at some given points. Let $G$ be an oriented graph with a length $\ell_e > 0$ associated to each edge $e$. Thus each edge $e$ may be seen as a segment of $\mathbb{R}$ and a differential operator $L_e = \frac{1}{2}a_e(x)\triangle + b_e(x)\nabla$ is associated to it. A differential operator $(L, \mathrm{Dom}(L))$ may be constructed from the $L_e$'s. Yet to properly define the operator $(L, \mathrm{Dom}(L))$, one needs to specify the behavior of the functions in the domain $\mathrm{Dom}(L)$ at the vertices. The domain $\mathrm{Dom}(L)$ of $L$ is composed of functions that are continuous on $G$, of class $\mathcal{C}^2$ on each edge and such that at any vertex $v$,

$$\sum_{e \sim v} \alpha_{v,e} \frac{\partial f}{\partial n_{v,e}}(v) = 0 \text{ with } \alpha_{v,e} \in [0,1], \ \sum_{e \sim v} \alpha_{v,e} = 1, \tag{57}$$

where $e \sim v$ means that $v$ is an endpoint of the edge $e$ and $\frac{\partial f}{\partial n_{v,e}}$ is the derivative of $f$ in the direction of the edge $e$ at the vertex $v$.

The condition (57) prescribes the flux in each direction, and is clearly a transmission condition. The one-dimensional problem (56) can be seen as a diffusion equation on a graph whose vertices are the $x_i$'s and whose edges are the intervals $(x_i, x_{i+1})$.

Indeed this generalization goes farther, since one may derive an Itô formula for this diffusion process with a local time at each vertex [33]. Besides, approximations by random walks is considered in [26]. The differential operator $(L, \mathrm{Dom}(L))$ is the infinitesimal generator of a stochastic process $X$, which can be seen as the limit of a diffusion process moving in small tubes around the edges and reflected at the boundary: See [35]. Diffusions on graphs are useful to describe the limiting behavior of perturbed Hamiltonian systems (see [34] and the subsequent works). Some of the Monte Carlo methods presented below in Section 11.9.3 may be adapted to this case. On some particular graphs, it is also possible to compute explicitly the transition density function of $X$: see [65] for example.

### 11.9.3. Numerical simulations

The question of the simulation of the solution of some SDE with discontinuous diffusion coefficients is a problem which have been hardly treated (See however [45, 16, 93]). Until recently, the case of processes generated by a divergence form operator of type $\frac{1}{2}\partial_{x_i}(a_{i,j}\partial_{x_j})$ have not been rigorously investigated, and the simulation of the process it generates in the multi-dimensional case remains an open problem (one knows from [80] that the diffusion may be approximated by a Markov chain, but the computation of the transition probabilities of this Markov chain is intractable in the general case). Due to many fields of applications (geophysics, electro/magneto-encephalography, ecology, astrophysics, ...), it is however of practical importance, and the main problem is related to the understanding of the behavior of the diffusion process when it reaches a hypersurface where $a$ is discontinuous.

In dimension one, Proposition 18 together with Proposition 17 allows us to describe exactly what happens to a diffusive particle that reaches a point where either $\rho$ or $a$ is discontinuous, if its trajectory is a realization of a trajectory of the process generated by the differential operator $L$ given by (54). Although one may think to use Corollary 2 to regularize the coefficients, this leads to unstable simulations (See [29] for numerical examples). Moreover, in practical situations, the jumps of the coefficients may be very high (for example, in a fissured porous media, the order of magnitude may be 1.000 or higher). In addition, it does not explain what happens to the diffusion.

In dimension one, several schemes to simulate diffusion processes with discontinuous coefficients may be given. and we present here briefly several approaches. In many cases, we have assumed that the coefficients $a$ and $\rho$ belong to BV whose points of discontinuity have no cluster points, or that they can be approximated well by such coefficients. We also note that one may assume that $b = 0$ by transforming the coefficients $a$ and $\rho$ into $a(x)\exp(-h(x))$ and $\rho(x)\exp(h(x))$ with $h(x) = 2\int_0^x b(z)/a(z)\rho(z)\,\mathrm{d}z$. This is the Zvonkin transform, as we saw it in the proof of Proposition 8.

We present several schemes. Once the behavior of the diffusion is understood, it becomes possible to provide other schemes or to mix them.

**A.** If there is one discontinuity, one may easily adapt the method proposed by E. Hausenblas in [42] for the reflected diffusion processes. There, the simulation skips the small excursions when the process reaches the boundary, and the particle jumps according to the entrance law of the diffusion process of excursions whose size is greater than a fixed parameter $\varepsilon$. The adaptation can be done using the decomposition of excursions measures given in [55], and by choosing the sign of the excursion whose size is greater than $\varepsilon$ with an independent Bernoulli random variable.

**B.** Using Proposition 17, one may find some one-to-one function $F$ to get rid of the term with the local time in (55). There, one obtains some SDE with a drift term and a diffusion term whose coefficients are discontinuous. It is then possible to apply results on the simulations of SDEs with discontinuous coefficients. This is the strategy chosen by M. Martinez and D. Talay [61, 62], where he uses the Euler scheme whose convergence is shown in [93]. In addition, the speed of convergence is computed.

**C.** Still using Proposition 17, after having approximated the coefficients by piecewise constant coefficients (after the drift have been removed), one may find a one-to-one function $G$ that transform the solution $X$ to (55) into the solution to the SDE

$$Y_t = G(X_0) + B_t + \sum_{n \in J} \frac{\sqrt{a(x^n+)/\rho(x^n+)} - \sqrt{a(x^n-)/\rho(x^n-)}}{\sqrt{a(x^n+)/\rho(x^n+)} + \sqrt{a(x^n-)/\rho(x^n-)}} L_t^{G(x^n)}(Y), \quad (58)$$

where $\{x^n\}_{n \in J}$ is the set of points where $a$ and $\rho$ are discontinuous. Hence, around each point $x^n$, $Y$ behaves like a Skew Brownian motion. The idea is then the following: We construct a sequence of points $\{y^k\}_{k \in K}$ such that $\{G(x^n)\}_{n \in J} \subset \{y^k\}_{k \in K}$ and $y_{k-1} < y_k < y_{k+1}$. We then set

$$\theta^i = \inf \left\{ t > \theta^{i-1} \, \middle| \, Y_t \in \{y^{k-1}, y^{k+1}\} \right\}$$

if $Y_{\theta^{i-1}} = y^k$. For this, we assume that $\theta^0 = 0$ and $Y_{\theta^0} = Y_0$ belong to $\{y^k\}_{k \in K}$. It is then possible to simulate exactly the Markov chain $(\theta^i, Y_{\theta^i})$. If $Y_{\theta^i} = y^k$ does not belong to $\{x^n\}_{n \in J}$, then it corresponds to simulate the first exit time and position from $[y^{k-1}, y^k]$ for the Brownian motion, which can be done numerically with the help of analytical expressions of the density of the Brownian motion killed when it exit from some interval. If $Y_{\theta^i} = y^k = x^n$ for $n \in J$, then we have to assume that the $\{y^k\}_{k \in K}$ are such that $x^n - y^{k-1} = y^{k+1} - x^n$ (this is why we use a finer grid $\{y^k\}_{k \in K}$). Then, using the "symmetry" of the negative and positive excursions of the Skew Brownian motion, the first exit time and position for $Y$ from $[y^{k-1}, y^{k+1}]$ can be deduced from the ones of the killed Brownian motion and the skewness parameter of $Y$ around $G(x^n)$.

Let us note that it is also possible to simulate $Y_T$ for an arbitrary time $T$. For this, we decide at each step with a Bernoulli random variable whether $\{\theta^{i-1} + \theta^i \leq T\}$ or not. In the former case, it is also possible to simulate $(\theta^i, Y_{\theta^i})$ given $\{\theta^{i-1} + \theta^i \leq T\}$ when one knows $(\theta^{i-1}, Y_{\theta^{i-1}})$. In the latter case, we are interested in simulating $Y_T$ given $\theta^i > T - \theta^{i-1}$, which can also be done from the density of the Brownian motion killed when it exit from some interval.

Once we have simulated the $(\theta^i, Y_{\theta^i})_{i=0,1,...}$, it is easy to use $G$ to get the distribution $(\tau \wedge T, X_{\tau \wedge T})$, there $\tau$ is the first exit time from some interval, and $T > 0$ is an arbitrary time.

This approach is presented in [58, 61].

**D.** We have already noted in Section 11.9.3 that J.-F. Le Gall has proved a Donsker theorem of SDEs of type $\mathrm{d}Z_t = \mathrm{d}B_t + \int_{\mathbb{R}} \nu(\mathrm{d}x) L_t^x(Z)$. This means that he constructed a random walk $(S_k)_{k \in \mathbb{N}}$ such that $n^{-1} S_{\lfloor n^2 t \rfloor}$ converges in distribution to $Z_t$ for any $t \geq 0$. The probability transitions of $S_k$ can be explicitly deduced from the measure $\mu$.

In [27], P. Étoré used the previous map $G$ to reduce the simulation of $X$ to this case, and then approximate the diffusion $Y$ using the random walk. He also computed the speed of convergence of this random walk. In some sense, this is also the simplification of the algorithm presented in **C.** where the $\{y^k\}_{k\in K}$ are equally spaced (at size $1/n$) and the exit times $\theta^{i+1} - \theta^i$ are replaced by their expectations, which are constant over all the position and equal to $1/n^2$.

**E.** In [28] (see also [29]), we constructed a bi-dimensional Markov chain $(\theta^k, Z^k)_{k\in\mathbb{N}}$ from a realization of a trajectory of $X$ generated by $L$ by setting

$$\tau^{k+1} = \inf\left\{\, t > \tau^k \,\middle|\, X_t \in \mathcal{G} \setminus \{X_{\tau^k}\} \,\right\}$$
$$\text{and } \theta^{k+1} = \theta^k + \mathbb{E}[\,\tau^{k+1} - \tau^k \mid (X_{\tau^k}, X_{\tau^{k+1}})\,], \; Z^{k+1} = X_{\tau^{k+1}},$$

where $\mathcal{G}$ is an arbitrary grid whose distance between two points is not necessarily equal to a fixed parameter.

Let us set

$$Z(t) = Z^k + \frac{t - \theta_k}{\theta_{k+1} - \theta_k}(Z^{k+1} - Z^k)$$

for $t \in [\theta^k, \theta^{k+1}]$ and $K(t) = \inf\left\{\, k \in \mathbb{N} \,\middle|\, \theta^k \geq t \,\right\}$. When the mesh of $\mathcal{G}$ decreases to 0, one may show that $(\theta^{K(t)}, Z(t))$ converges uniformly in $t \in [0, T]$ in probability to $(t, X_t)$ and the rate of convergence may be computed. This construction no longer uses the SDE (55), and the assumption that $a$ and $\rho$ belong to BV may be dropped.

We then approximate $X$ by simulating a Markov chain with the distribution of $(Z^k, \theta^k)_{k=0,1,2,\dots}$. For this, the probability transition of $(Z^k)_{k\in\mathbb{N}}$ is deduced from the scale function, since for

$$v_i(x) = \frac{S(x_{i-1}) - S(x)}{S(x_{i+1}) - S(x_{i-1})} \text{ for } x \in (x_{i-1}, x_{i+1})$$

where $x_{i-1} < x_i < x_{i+1}$ are three successive points on the grid $\mathcal{G}$, then $\mathbb{P}[\, Z^{k+1} = x_{i+1} \mid Z^k = x_i\,] = v_i(Z^k)$. Let us note that $v_i(x)$ is solution to $Lv_i(x) = 0$ on $(x_{i-1}, x_{i+1})$ with $v_i(x_{i-1}) = 1$ and $v_i(x_{i+1}) = 1$. The value of $\theta^k$ is then computed by solving $Lu_i(x) = -\widetilde{v}_i(x)$ on $(x_{i-1}, x_{i+1})$ with $u_i(x_{i-1}) = u_i(x_{i+1}) = 0$ and $\widetilde{v}_i = v_i$ if $Z^{k+1} = x_{i+1}$ and $\widetilde{v}_i = 1 - v_i$ if $Z^{k+1} = x_{i-1}$. In facts, $\theta^k = v_i(Z^k)$.

**F.** In [19, 22], M. Decamps, A. De Schepper, M. Goovaerts and W. Shouten compute the densities of some process with discontinuous coefficients, in view of financial applications. For this, the Fokker-Planck equation is solved using perturbation formula of the density with respect to the density of the SBM: For the solution to the SDE

$$\mathrm{d}X_t = \sigma(X_t)\,\mathrm{d}B_t + b(X_t)\,\mathrm{d}t + (2\alpha - 1)\,\mathrm{d}L_t^{x^*}(X),$$

where $B$ is a Brownian motion, the write the density $p(t, x, y)$ of $X$ as

$$p(t, x, y) = \frac{1}{\sigma(x)} \exp\left(\int_{y_0}^{y} b(z)\,\mathrm{d}z\right) p^\beta(t, y_0, y)$$
$$\times \mathbb{E}_{y_0}\left[\exp\left(-\int_0^t c(X_s^\beta)\,\mathrm{d}s - \gamma L_t^{\Phi(x^*)}\right)\middle| X_t^\beta = y\right], \quad (59)$$

where $X^\beta$ is the SBM($\beta$) with

$$\beta = \frac{(\alpha + 1/2)\sigma(x_-^\star) + (\alpha - 1/2)\sigma(x_+^\star)}{\sigma(x_-^\star) + \sigma(x_+^\star)},$$

$p^\beta(t, x, y)$ is its density,

$$\Phi(x) = \int^x \frac{\mathrm{d}z}{\sigma(z)} \text{ and } c(y) = \frac{b(\Phi^{-1}(y))}{\sigma(\Phi^{-1}(y))} - \frac{1}{2}\sigma'(\Phi^{-1}(y)).$$

They then consider several "Skew models" — such as the Self Exciting Threshold (SET) Cox-Ingersoll-Ross, SET Vasicek, SET Libor market, ... — and use spectral decompositions to compute and/or approximate the density given in (59).

They also consider the case of Skew Bessel processes.

### 11.9.4. Parameter estimation

In [1, 61], O. Bardou and M. Martinez have constructed a scheme to estimate both the skewness parameter $\gamma$ and the position of the point $x$ of the doubly reflected Skew Brownian motion $\mathrm{d}X_t = \mathrm{d}B_t + \gamma \mathrm{d}L_t^x(X) + \mathrm{d}L_t^{-1}(X) - \mathrm{d}L_t^1(X)$. Their construction relies on the ergodicity of the underlying process, and its connections with PDEs with discontinuous coefficients.

## 11.10. Multi-dimensional extensions

The natural multi-dimensional extension of the Skew Brownian motion is that of a SDE involving the local time of the process in some hyper-surface. Yet it is difficult task to construct such a process. We then present briefly the method proposed by N. Portenko on generalized diffusion processes, where the semigroup of the process is constructed using a Volterra series. Afterwards, we give some references to other constructions, some of them relying on fixed point theorems on SDEs, and other on Dirichlet forms.

### 11.10.1. Generalized diffusion processes

In [69, 70, 71], N. Portenko develops the concept of diffusions with *singular drift*, which we already used in Section 2.

Let $S$ be a surface in $\mathbb{R}^d$ which separates $\mathbb{R}^d$ into two regions: the interior region $D$ and the exterior region $\mathbb{R}^d \setminus D$, $N$ a vector field on $S$, $q$ a continuous function from $S$ to $\mathbb{R}$. We assume that $S$ satisfies both the interior and the exterior sphere property, which means that for any $x \in S$, one can find some balls $B$ and $B'$ with $B \subset \overline{D}$, $B' \subset \overline{\mathbb{R}^d \setminus D}$, $B \cap S = B' \cap S = \{x\}$.

Let also $a$ be a continuous function with values in the space of $d \times d$-symmetric matrices, which is uniformly elliptic and bounded. A *generalized diffusion* is a diffusion process whose infinitesimal generator $L$ is formally written

$$L = \frac{1}{2} \sum_{i,j=1}^{d} a_{i,j}(x) \frac{\partial^2}{\partial x_i \partial x_j} + \sum_{i=1}^{d} \delta_S(x) q(x) N_i(x) \frac{\partial}{\partial x_i}. \tag{60}$$

The idea is then to construct first the semi-group $(P_t)_{t>0}$ to $L$, it to show it generates a Feller process $(X, \mathbb{P}_x)$. The solution of the parabolic PDE

$$\frac{\partial u(t,x)}{\partial t} = Lu(t,x) \text{ on } \mathbb{R}_+^* \times \mathbb{R}^d \text{ with } u(0,x) = \varphi(x) \tag{61}$$

is then $u(t,x) = \int_{\mathbb{R}^d} p(t,x,y)\varphi(y)\,\mathrm{d}y = \mathbb{E}_y[\varphi(X_t)]$.

The problem is of course to interpret the meaning of the $\delta$-function in the drift. One may see (61) as a PDE with a transmission condition on $S$, as (61) may be written

$$\begin{cases} \frac{\partial u(t,x)}{\partial t} = Lu(t,x) \text{ on } \mathbb{R}^d \setminus S, \\ (1+q(x))N_+(x) \cdot \nabla u(t,x) = (1-q(x))N_-(x) \cdot \nabla u(t,x) \\ u(t,\cdot) \text{ is continuous on } S, u(0,\cdot) = \varphi \text{ on } \mathbb{R}^d, \end{cases} \tag{62}$$

where $N_+$ and $N_-$ are the inner and outer conormal to $S$, that is $N_+(x) = a(x)n(x)$ and $N_-(x) = -N_+(x)$, where $n(x)$ is the normal unit vector to $S$ at $x$, which is directed inward. There, $N_\pm(x) \cdot \nabla u(t,x) = \lim_{\varepsilon \to 0} \varepsilon^{-1}(u(t,x \pm \epsilon N_\pm(x)) - u(t,x))$.

As the $\delta$-function imposes some flux condition on $S$, we can think to use a *single layer potential*, which is generated by a distribution of charges $V(t,\cdot)$ in $S$ at time $t$ which we choose to be

$$V(t,x) = \frac{q(x)}{2}\left(N_+(x) \cdot \nabla u(t,x+) + N_-(x) \cdot \nabla u(t,x-)\right), \ x \in S. \tag{63}$$

Let $p^0(t,x,y)$ be the fundamental solution to $L^0 = \frac{1}{2}a_{i,j}\partial^2_{x_i x_j}$. The potential generated by $V(t,x)$ on $\mathbb{R}^d$ is then equal to

$$v(t,x) = \int_0^t \int_S p^0(t-s,x,y)V(s,y)\,\mathrm{d}\sigma_y\,\mathrm{d}s,$$

where $\sigma_y$ is the Lebesgue measure on $S$, and according to classical results on potential theory, $v$ is harmonic on $\mathbb{R}^d \setminus S$ and continuous on $\mathbb{R}^d$ (and then when crossing $S$). In addition, for $x \in S$,

$$N_\pm(x) \cdot \nabla v(t,x) = \int_0^t \int_S N_+(x) \cdot \nabla p^0(t-s,x,y)V(s,y)\,\mathrm{d}\sigma_y\,\mathrm{d}s \mp V(t,x). \tag{64}$$

As $v(0,x) = 0$, we then sought $u$ to be the superposition of $v$ and $u^0(t,x)$ defined by

$$u^0(t,x) = \int_{\mathbb{R}^d} p^0(t,x,y)\varphi(y)\,\mathrm{d}\sigma_y. \tag{65}$$

This solution corresponds to the case $q = 0$ and satisfies on $S$, $N_+(x)\nabla u^0(t, x+) = N_-(x)\nabla u^0(t, x-)$.

By combining the previous equation with Equations (63), (64), we obtain that in order that the flux condition in (62) is satisfied, the charge distribution $V$ shall be solution to

$$V(t, x) = \int_{\mathbb{R}^d} N_+(x) \cdot \nabla p^0(t, x, y)\varphi(y) \, \mathrm{d}y$$
$$+ \int_0^t \int_S N_+(x) \cdot \nabla p^0(t, x, y)V(s, y)q(y) \, \mathrm{d}\sigma_y \, \mathrm{d}s. \quad (66)$$

This equation may be solved using Volterra series, which means by constructing recursively

$$V^{(0)}(t, x) = \int_{\mathbb{R}^d} N_+(x) \cdot \nabla p^0(t, x, y)\varphi(y) \, \mathrm{d}y,$$

$$V^{(n+1)}(t, x) = \int_0^t \int_S N_+(x) \cdot \nabla p^0(t, x, y)V^{(n)}(s, y)q(y) \, \mathrm{d}\sigma_y \, \mathrm{d}s$$

and establishing the convergence of the series $V(t, x) = \sum_{n \geq 0} V^{(n)}(t, x)$ (the advantage we get in Section 2 was that there was no need to consider such a series, since $N_+(x) \cdot \nabla p^0(t, 0, y) = 0$ when $p^0(t, x, y)$ is the one-dimensional heat kernel). We have then constructed a charge distribution on $S$ whose single layer potential $u$ gives a solution to (61) outside $S$. It is then easily obtained that $u$ is solution to (62).

We are now willing to construct the semi-group $(P_t)_{t>0}$ of $L$. For this, we set $P_t^0\varphi(x) = \int_{\mathbb{R}^d} p^0(t, x, y)\varphi(y) \, \mathrm{d}y$ and

$$P_t\varphi(x) = P_t^0\varphi(x) + \int_0^t \int_S p^0(t - s, x, y)V(s, y)q(y) \, \mathrm{d}\sigma_y \, \mathrm{d}s,$$

where $x \in \mathbb{R}^d$, $t > 0$ and $V$ is the solution to (66). This is a perturbation formula. It is then possible to prove that $P_{s+t}\varphi(x) = P_s P_t \varphi(x)$ for any $x \in \mathbb{R}^d$, that $P_t\varphi(x) \xrightarrow[t \to 0]{} \varphi(x)$ for $x \in \mathbb{R}^d$ and that

$$N_+(x) \cdot \nabla P_t\varphi(x) = (1 - q(x))V(t, x) \text{ and } N_-(x) \cdot \nabla P_t\varphi(x) = (1 + q(x))V(t, x),$$

so that $(t, x) \mapsto P_t\varphi(x)$ is solution to (62).

As for the proof of Proposition 1, if $|q(x)| \leq 1$, then $P_t\varphi$ is non-negative when $\varphi$ is non-negative. In addition $P_t 1 = 1$, so that $(P_t)_{t>0}$ is a contraction semi-group for the uniform norm.

**Theorem 15** ([69, Theorem 1])**.** *Under the previous hypotheses on $S$, $a$ and $q$ and if $|q(x)| \leq 1$ for all $x \in S$, The semi-group $(P_t)_{t>0}$ generates a conservative, continuous Markov process $(X, \mathbb{P}_x, x \in \mathbb{R}^d)$.*

*Remark* 17. If $S$ is some hyper-plane in $\mathbb{R}^n$ and $a$ is the identity matrix, then similar computations also lead to the same result [69, Section 7].

One may wish to characterize the diffusion process constructed so by its diffusion and drift coefficients. It is standard that the diffusion coefficient is characterized by $a(x)\theta \cdot \theta = \lim_{t \to 0} \mathbb{E}_x[((X_t - x) \cdot \theta)^2]$ for any $\theta \in \mathbb{R}^d$. Using the semi-group and a continuous test function $\psi$ with compact support, one can show that

$$\int_{\mathbb{R}^d} \psi(x) \frac{1}{t} \int_{\mathbb{R}^d} ((y - x) \cdot \theta)^2 P(t, x, \mathrm{d}y) \, \mathrm{d}x \xrightarrow[t \to 0]{} \int_{\mathbb{R}^d} \psi(x) a(x) \theta \cdot \theta \, \mathrm{d}x,$$

where $P(t, x, \mathrm{d}y)$ is the probability transition of $P$. Still using the test function $\psi$, one gets that is

$$\int_{\mathbb{R}^d} \psi(x) \frac{1}{t} \int_{\mathbb{R}^d} ((y - x) \cdot \theta) P(t, x, \mathrm{d}y) \, \mathrm{d}x \xrightarrow[t \to 0]{} \int_S \psi(x) q(x) (N_+(x) \cdot \theta) \, \mathrm{d}\sigma_x,$$

which characterized a drift concentrated on $S$.

In [70], N. Portenko decomposes the process $X$ as a semi-martingale.

**Theorem 16** ([70, Theorems 1 and 2])**.** *There exists a continuous additive functional* $(\zeta_t)_{t \geq 0}$ *of bounded variation such that* $M_t = X_t - X_0 - \zeta_t$ *is a continuous square integrable martingale with* $\langle M \rangle_t = \int_0^t a(X_s) \, \mathrm{d}s$ *and such that* $\frac{\mathrm{d}}{\mathrm{d}t} \mathbb{E}_x[\zeta_t^i] = T_t(q(x) N_+^i(x))$ *and* $\zeta_t^i = \int_0^t \mathbf{1}_{\{X_s \in S\}} \, \mathrm{d}\zeta_s^i$ *for* $i = 1, \ldots, d$.

Some of these results have been extended recently (see for example [72]), for example to take into account some elastic killing effect when the process passes through the surface $S$.

### 11.10.2. SDE with local time in the multi-dimensional case

There exists several works on SDEs with local time that may appears as generalizations of the SBM. However, their constructions are rather technical in general, and the hypotheses tedious to write. This is why we give here just a few references.

**A.** The article [84] provides a construction of some diffusion $X$ associated to the Dirichlet form

$$(u, v) \in \mathrm{H}_0^1(\mathbb{R}^d) \mapsto \int_{\mathbb{R}^d} a(x) \nabla u(x) \nabla v(x) (\alpha \mathbf{1}_G(x) + (1 - \alpha) \mathbf{1}_{\mathbb{R}^d \setminus G}(x)) \rho(x) \, \mathrm{d}x$$

where $G$ is a domain of $\mathbb{R}^d$, $\alpha \in (0, 1)$, $a$ is uniformly elliptic, $\rho > 0$ and of course, $a$ and $\rho$ have sufficient integrability condition. If $a$ and $\rho$ are "weakly differentiable", then one can decompose $X$ as some SDE with a local time associated to $\partial G$.

**B.** Given a measure $\mu$ which is singular with respect to the Lebesgue measure, Y. Ōshima proved in [66] the weak existence and uniqueness of the $d$-dimensional process solution to

$$\mathrm{d}X_t^i = \sum_{j=1}^d \left( \sigma_{i,j}(X_t) \, \mathrm{d}B_t^j + b_i(X_t) \, \mathrm{d}t + \tau_{i,j}(X_t) \, \mathrm{d}M_t^j + \beta_i(X_t) \ell_t^\mu \right), \; i = 1, \ldots, d,$$

where $\sigma\sigma^{\mathrm{T}} = (a_{i,j})_{i,j=1}^d$, $\tau\tau^{\mathrm{T}} = \alpha/a_{1,1}$, $a_{1,1} \geq c > 0$, $b_i = \frac{1}{2}\sum_{j=1}^d \partial_{x_j} a_{j,i}$, $\beta_i = \frac{1}{2a_{1,1}}\sum_{j=2}^d \partial_{x_j}\alpha_{j,i}$, $\beta_1 = \tau_{1,j} = \tau_{i,1} = 0$, $B$ is a $d$-dimensional Brownian motion, $\ell_t^\mu = \int_{\mathbb{R}} \ell_t^x \mu(\,\mathrm{d}x_1)$, where $(\ell^x)_{x\in\mathbb{R}}$ is a family of continuous non-decreasing processes such that $\ell_t^x = \int_0^t \mathbf{1}_{\{X_s^1=x\}}\,\mathrm{d}\ell_s^x$ and

$$\int_{\mathbb{R}} \ell_t^x f(x)\,\mathrm{d}x_1 = \int_0^t f(X_s^1) a_{1,1}(X_s)\,\mathrm{d}s$$

for all continuous function $f$ with compact support, $M$ is a family of continuous local martingales with $\langle M^i, M^j\rangle_t = \delta_{i,j}\ell_t^\mu$ and $\langle M^i, B^j\rangle = 0$ for $i,j = 1,\dots,d$. The coefficients $\sigma$ and $\tau$ are assumed to be Lipschitz continuous.

Y. Ōshima also proves that this solution $X$ is generated by the Dirichlet form

$$(u,v) \mapsto \frac{1}{2}\int_{\mathbb{R}^d} a(x)\nabla u(x)\nabla v(x)\,\mathrm{d}x + \frac{1}{2}\int_{\mathbb{R}^d} \alpha_{i,j}(x)\nabla u(x)\nabla v(x)\eta(\,\mathrm{d}x),$$

where $\eta$ is the measure $\eta(\,\mathrm{d}x) = \mu(\,\mathrm{d}x_1)\,\mathrm{d}x_2\cdots\mathrm{d}x_d$, which is then the superposition of two Dirichlet forms[2].

**C.** The previous construction was extended by S. Takanobu in [81, 82], where the assumptions that $\beta_1 = 0$ and $a_{1,1} \geq c > 0$ are replaced by the assumption that $\sup_{y\in\mathbb{R}^{d-1}} |\beta_d(y,x)\mu(\{x\})| < 1$ for all $x \in \mathbb{R}$ and $\sum_{j=1}^d \sigma_{j,d}^2 \geq c > 0$ for some positive constant.

**D.** In [77], A.-S. Sznitman and S.R.S. Varadhan constructs the solution to the SDE

$$X_t = x + B_t + \sum_{k=0}^N V_k L_t^0(n_k \cdot X),$$

where the $n_k$'s are unit vectors generating distinct hyperplanes $H_k$, and $V_k$ are constant vectors with $V_k \cdot n_k = 0$. This appears as a way to construct the solution of some $N$ dimensional diffusion corresponding to $N$ one-dimensional particles in interaction.

**E.** In [95, 96, 97], L. Zaitseva shows the strong existence and uniqueness of the SDE

$$X_t = x + B_t + \int_0^t \sigma(\pi_S(X_s))\,\mathrm{d}\widetilde{B}_{L_s} + \int_0^t (\beta\nu + \alpha(\pi_S(X_s)))\,\mathrm{d}L_s, \tag{67}$$

where $S$ is a hyperplane of $\mathbb{R}^d$ orthogonal to the unit vector $\nu$, $\pi_S$ the projection from $\mathbb{R}^d$ onto $S$, $B$ a Brownian motion on $\mathbb{R}^d$, $\widetilde{B}$ a Brownian motion on $S$, $\beta \in (-1,1)$, $(L_s)_{s\geq 0}$ is the local time of $\nu \cdot X$, $\alpha$ a measurable function from $S$ to $S$ and $\sigma$ a measurable function such that $\sigma\sigma^{\mathrm{T}}$ takes its values in the space of symmetric matrices. The functions $\sigma$ and $\alpha$ are assumed to be bounded and

---

[2]The article [83] and subsequent works also contains some constructions of diffusion processes associated to the superposition of Dirichlet forms defined on hypersurfaces, but without giving some SDE solved by them.

Lipschitz continuous. This diffusion is shown in [97] to be a generalized diffusion in the sense of N. Portenko [71] with drift $(\beta\nu + \alpha \circ \pi_S)\delta_S$ and a diffusion term $\mathrm{Id} + \sigma \circ \pi_S \delta_S$.

### 11.11. Domains of applications

Being related to diffusion processes with discontinuous coefficients and to diffusions on graphs, the Skew Brownian motion has potential applications in many fields of physics: in astrophysics [98], in geophysics and study of heterogeneous media [35, 50, 56, 57, 73, 92], in perturbed Hamiltonian systems [37], study of hysteresis phenomena [32], in biology (some modelling of nerve impulses involves parabolic PDEs as in Section 11.9.2 [64]), ecology (the dynamic of populations moving between different refuges in studied in [15]), in finance and actuarial sciences [19, 20, 21], in optimization (see below)...

In all these applications, the properties of the Skew Brownian motion may be used in order to solves the related problems by Monte Carlo methods, or to use its properties to solve such a problem analytically.

#### Applications to optimal problems

The article [53] looks for a a drift $b \geq 0$ with $\int_0^1 b(x)\,\mathrm{d}x < +\infty$ that minimizes $\mathbb{E}_{x_0}\tau$, where $\tau = \inf\left\{ t \geq 0 \,\middle|\, X_\tau = 1 \right\}$ and $X$ is the solution to $\mathrm{d}X_t = \mathrm{d}B_t + b(X_t)\,\mathrm{d}t$ with a reflecting boundary at 0. Indeed, the optimal drift is a measure that can be computed explicitly, and that is of type $b(\mathrm{d}x) = \widetilde{b}(x)\mathrm{d}x + \kappa\delta_{x_0}(\mathrm{d}x)$. Hence it involves a process of Skew Brownian motion type.

The constructions in [2] and in [17] are also related to methods for solving resources' allocations problems (see [60] for example). Recently, the variably skew Brownian motion presented in Section 11.6 has also been used for solving the Skorokhod embedding problem [18].

**Acknowledgment.** The author is grateful to Pierre Étoré for interesting discussions about the Skew Brownian motion, especially about the approach by N. Portenko. The author also wishes to thank the referees who helped him to improve the quality of this article.